\input amstex.tex
\documentstyle{amsppt}      
\baselineskip=10pt
\parindent=0.2truein
\pagewidth{6.6truein}
\pageheight{9.0truein}
\topmatter
\endsubjclass

\title  Remarks on Ruelle Operator and\\
invariant line fields problem
\endtitle
\author Peter M. Makienko
\endauthor
\address{Permanent address: Institute for Applied Mathematics, \newline  9 Shevchenko str.,
Khabarovsk, Russia and \newline Instituto de Matematicas,\newline 
Av. de Universidad s/N., Col. Lomas de Chamilpa, C.P. 62210, Cuernavaca, Morelos, Mexico }
\endaddress
\thanks{This work has been partially supported by the Russian Fund of Basic
Researches, Grant 99-01-01006.}
\endthanks
\endtopmatter

\document

\heading{\bf Introduction and main statements}\endheading
\par Let $ R  $ be a rational map with non-empty Fatou set. Assume that $ J(R) $ supports an invariant non trivial conformal structure 
$ \mu.$ Let $ f_\mu $ be its respective quasiconformal map. The main idea of this work is to find the conditions which allow to construct 
a quasiconformal map $ h $ supported already by the Fatou set, so that $ h $ and $ f_\mu $ generate the same infinitesimal deformation of
 $ R $ (see also {\cite {\bf Mak}}). This approach allows us to formulate conditions (in terms of Ruelle-Poincare series) which guarantee
 the absence of non trivial invariant conformal structures on the Julia set, see the theorems below. Also the necessary and sufficient 
conditions (in terms of convergence of sequences of measures) of existence of invariant conformal structures on $ J(R) $ are obtained.
\subheading{Main results}
\par  Let $ R $ be a rational map with non-empty Fatou set $ F(R).$
Denote by $ Pc(R) $ the postcritical set of $ R. $ In further, we always 
 suppose that
 \roster
 \item all critical points are simple (that is if $ R'(c) = 0, $ then
 $ R''(c) \neq 0 $),
 \item there exist no simple critical relations (that is $ R $ has a simple critical relation 
 iff there exist integers $ n $ and $ m $ and two {\it different critical points}
  $ c_1 $ and $ c_2, $ so that the following equality
 $$ R^n(c_1) = R^m(c_2) $$
 is hold)
 \endroster

\par  Start again with a rational map $ R $ and consider two actions $ R^*_{n,m} $ and $ {R_*}_{n,m} $ on a function $ \phi $ at point 
$ z $ by the formulas
$$
\split
 &R^*_{n,m}(\phi) = \sum\phi(J_i)(J_i')^n(\overline{J_i'})^m = 
\sum_{y \in R^{-1}(z)}\frac{\phi(y)}{(R'(y))^n(\overline{R'(y)})^m}, \text { and }\\
&{R_*}_{n.m}(\phi) = \phi(R)\cdot (R')^n\cdot\overline{(R')^m},
\endsplit
$$
where $ n $ and $ m $ are integers and $ J_i, i = 1, ..., d $ are branches
of the inverse map $ R^{-1}. $ Then we have 
$$ R^*_{n,m}\circ {R_*}_{n,m}(\phi) = deg(R)\cdot\phi.
$$
In other words the actions above  are the natural action 
of $ R $ on the spaces of the forms of the kind $ \phi(z)Dz^mD\overline{z}^n. $ 
\proclaim{Definition}
\roster
\item The operator $ R^* = R^*_{2,0} $ is called {\rm Ruelle operator} of rational map $R.$
\item  The operator 
$ \vert R^*\vert = R^*_{1, 1} $ is called {\rm modulus of the Ruelle operator}
\item The operator  $ B_R = {R_*}_{-1,1} $ is called {\rm Beltrami operator} of rational map $ R. $
\endroster
\endproclaim
 Operators $ R^* $ and $ \vert R^*\vert $ and its right inverse $ R_* = \frac{{R_*}_{2,0}}{deg(R)} $ and $ \vert R_*\vert(\phi) = \frac{\vert {R_*}_{2,0}\vert(\phi)}{deg(R)} = \frac{\phi(R)\vert R'\vert^2}{deg(R)} $ 
map the space $ L_1(\overline{\Bbb C}) $ into itself 
with the unite norm. The operator $ B_R $ maps the space $ L_\infty(\overline{\Bbb C}) $ into itself evidently with the unite norm.
\proclaim{Definition} Let $ R \in  {\Bbb C}P^{2d+1} $ be a rational map.
The component of {\bf J}-stability of $ R $ is the following space.
$$
\aligned qc_J(R) = \bigl\{
&F \in {\Bbb C}P^{2d+1}: \text{ there are neighborhoods } U_R
\text{ and } U_{F} \text{ of }{ \bold J}( R)\,\, \text{ and }
{\bold J}(F),\\
&\text{respectively and a quasiconformal homeomorphism } h_F : U_R \rightarrow
 U_{F} \text{ such that }\\
&F = h_F \circ R \circ {h^{-1}_F}\bigr\}\big/PSL_2(\Bbb C).
\endaligned 
$$
\endproclaim
\proclaim{Definition} The space of invariant conformal structures or Teichmuller space $ T(J(R)) $on $ J(R) $ for a given rational map $ R $ is the following space
$$ 
T(J(R)) = \{Fix(B_R)(L_\infty(J(R))\}.
$$
where $ Fix(B_R)(L_\infty(J(R)) $ is the space of fixed elements of Beltrami operator $ B_R: L_\infty(J(R))\rightarrow L_\infty(J(R)).$
Due to D. Sullivan (see \cite{\bf S}) the dimension of $ T(J(R)) $ is bounded above by $ 2deg(R) - 2.$
\endproclaim
\proclaim{Definition} We will call a rational map $ K(z) $ {\rm convergent} iff there exists a rational map $ R \in qc_J(K) $ such that:
\roster
\item there exists a point $ a \in \overline{\Bbb C} $ with long orbit $ (\geq 2deg(R)) $,
\item for any $ x \in \cup_nR^n(a) $ there exists a sequence of integers $ \{N_i\} $ such that the expression
$$ A_{N_i}(x, R,) = \frac{1}{N_i}\sum_{j = 0}^{N_i - 1}\sum_k\frac{\vert (R^*)^j(\tau_x)(c_k)\vert}{\vert R''(c_k)(R^{N_i - j - 1})'(R(c_k))\vert}\tag{*}
$$
is uniformly bounded independently on $ i $ and where $ c_k $ are critical points for  map $ R.$
\endroster
\endproclaim 
\proclaim{Theorem A} Let $ R $ be a convergent map with simple critical points  and no simple critical relations. Assume that  Fatou set is non-empty and Lebesgue measure of postcritical set is zero. Then there is no non-trivial invariant conformal structures on $ J(R).$
\endproclaim
\proclaim{Definition} Ruelle-Poincare series.
\roster
\item {\rm Backward  Ruelle-Poincare series.}
$$ RS(x, R, a) = \sum_{n = 0}^\infty(R^*)^n(\tau_a)(x), $$
where  $ \tau_a(z) = \frac{1}{z - a} $ and $ a \in \overline{\Bbb C} $ is a parameter. The series  
$$ S(x, R) = \sum_{n = 0}^\infty\vert R^*\vert^n({\bold 1}_{\Bbb C})(x)
$$
is called {\rm Backward Poincare series}.
\item {\rm Forward Ruelle-Poincare series.}
$$ RP(x,R) = \sum_{n = 0}^\infty \frac{1}{(R^n)'(R(x))}.
$$
The series
$$ P(x,R) = \sum_{n = 0}^\infty \frac{1}{\vert(R^n)'(R(x))\vert}
$$
is called {\rm forward Poincare series.}
The series
$$ A(x, R, a) = \sum_{n = 0}\frac{1}{(R^n)'(a)(x - R^n(a))}
$$
is called {\rm modified Ruelle-Poincare series.}
\endroster
\endproclaim
Note that the Ruelle-Poincare series are a kind of generalizations of the Poincare series introduced by C. McMullen for rational maps (see \cite{\bf MM}). 

\proclaim{Corollary A} Let $ R $ be a rational map with simple critical points and no simple critical relation. Then $ R $ is convergent map if for some $ a \in {\Bbb C} $ with $\#\{\cup_iR^i(a)\} > 2deg(R) - 1 $ the one of the following is true.
\roster
\item Collet-Eckmann case. For any critical point $ c $ and an $ x \in \{\cup_iR^i(a)\}  $ the series $ RS(c,R, x) $ and $ RP(c, R) $ are absolutely convergent.
\item For any critical point $ c  $ and $ x \in \{\cup_iR^i(a)\}  $ one of the series  $ RS(c,R, x) $ or $ RP(c, R) $ is  absolutely convergent and the second one has uniformly bounded elements.
\item Conjectural case. Both series diverge slow enough (like harmonic series).
\endroster
\endproclaim
 As it will be shown below,  our definition of Collet-Eckmann maps (see item 1 above) is a generalization of one given by Feliks Przytycki (see {\cite{\bf P}}). Hence the item 1 of corollary A is reproof of Przytycki result in a weaker sense (we can not show that  $ m(J(R)) =0 $).
\par The third case is conjectural in the sense that the case  looks like a $Log$-Collet-Eckmann case (that is $ \vert (R^n)'(R(c)\vert \sim C\cdot n$). It is not clear if it does exist a map with such behavior of its Ruelle-Poincare series.
\par The next proposition gives some formal relations between Ruelle-Poincare series. 
\proclaim{Definition}  We denote the Cauchy product of series $ A $ and $ B $ by $ A\otimes B $. Let us recall that if $ A = \sum_{i = 1} a_i $ and $ B = \sum_{i = 1} b_i, $ then $ C =  A\otimes B  = \sum_{i = 1} c_i, $ where $ c_i = \sum_{j = 1}^i a_jb_{i - j}.$ 
\endproclaim

Then we have the following proposition. 
\proclaim{Proposition A} Let $ R $ be a rational map with simple critical points. Let $\infty $ be a fixed point for $ R. $ Then there exist the following formal relations. 
$$ \align 
RP(a, R) - 1
&= \sum_i\lambda^i -\sum_i\frac{1}{R''(c_i)}RS(c_i,R,a)\otimes RP(c_i,R),\text{ where } \lambda \text{ is the multiplier of } \infty\\
RS(x, R, a) 
&=A(x, R, a) + \sum_k\frac{1}{R''(c_k)}A(c_k, R, a)\otimes  RS(x,R,R(c_k)),
\endalign
$$
where $ c_k $ are critical points of $ R. $
\endproclaim
   \par For polynomials of degree two this approach gives the following statement.
\proclaim{Theorem B} Let $ R(z) = z^2 +c $ and $ S_L = \sum^{L}_{j = 0}\frac{1}{(R^{j})'(c)}. $ Assume that there exists a subsequence $\{n_i\} $ of integers such that
\roster
\item $ \lim_{i \rightarrow \infty} \vert (R^{n_i})'(c)\vert =\infty $ and
 $\overline{\lim}_{i \rightarrow \infty} \vert S_{n_i}\vert > 0 $ or
\item $ \vert (R^{n_i})'(c)\vert \sim C=Const $ for  $ i\to\infty $ and 
$\overline{\lim}_{i \rightarrow \infty} \vert S_{n_i}\vert = \infty $ 
\endroster
Then there exists no invariant conformal structures on its Julia set.
\endproclaim

\par Finally, we give necessary and sufficient conditions for the existence of measurable invariant conformal structures on its Julia set in the case when the postcritical set has Lebesgue measure zero.
Let $ U $ be a neighborhood of $ J(R). $
We call $ U $ -{\it essential neighborhood} iff
\roster
\item $ U $ does not contain disks centered at all attractive and superattractive  points and
\item $ R^{-1}(U) \subset U. $
\endroster

\proclaim{Definition} Let us define the space $ H(U) \subset C({\overline U}), $ where $ C({\overline U}) $ is space of continuous functions and 
\roster
\item $ U $ is an essential neighborhood of $ J(R) $ and
\item $ H(U) $ consists of $ h\in C({\overline U}) $ such that $ \frac{\partial h}{\partial \overline{z}} $ (in the sense of distributions) belongs to $ L_\infty(U) $
\item $ H(U) $ inherits the topology of $ C({\overline U}).$
\endroster
\medskip
\par Measures $ \nu^i_l. $
\roster
\item Let $ c_i $ and $ d_i $ be critical points and critical values, respectively.
Then  define $ \mu_n^i = \frac{\partial}{\partial\overline{z}}((R^*)^n(\frac{1}{z - d_i})) $ (in sense of distributions). We will show below that $  (R^*)^n(\frac{1}{z - d_i})$ $ = $ $\sum_{j = 0}^n\frac{\alpha^i_j}{z - R^j(d_i)} $ and hence $ \mu_n^i = \sum_{j = 0}^n\alpha^i_j\delta_{R^j(d_i)}, $ where $ \delta_a $ denotes the delta measure with mass at the point $ a. $
\item Define by $ \nu^i_l $ the average $ \frac{1}{l}\sum_{k = 0}^{l - 1}\mu_k^i. $
\endroster
\endproclaim
\proclaim{Remark 1} For $ R(z) = z^2 + c $ we have 
$$ \mu_0 = \delta_c, \mu_i = \frac{\delta_{R(c)} - \delta_{c}}{R\prime(c)},$$
$$\mu_2 = \frac{\delta_{R^2(c)}}{(R^2)\prime(c)} - \frac{\delta_{c}}{(R^2)'(c)} - \frac{\mu_1}{R'(c)}, $$
$$ \mu_n = \frac{\delta_{R^n(c)}}{(R^n)'(c)} - \frac{\delta_{c}}{(R^n)'(c)} - 
\frac{\mu_1}{(R^{n - 1})'(c)} - ...- \frac{\mu_{n - 1}}{R'(c)}. $$
Now define by $ A $ the series $\sum_i \mu_i, $ then by above we have the following formal equality
$$ A\otimes RP(0, R) = \sum_{i =0}^{\infty} \frac{\delta_{R^i(c)}}{(R^i)'(c)}
$$
\endproclaim
\par In general the coefficients $ \alpha_j^i $ in definition above can be expressed as a combinations of the elements of Cauchy product $ RP(c_i,R,d_i)\otimes RP(c_i,R) $ of Ruelle-Poincare series.
\proclaim{Theorem C} Let $ R $ be a rational map with simple critical points and no simple critical relations. Assume that $ F(R) \neq \emptyset $ and $ m(Pc(R)) = 0, $ where $ F(R)$ is the Fatou set and $ m $ is denote the Lebesgue measure. Then $ T(J(R)) = \emptyset $ if and only if there exist an essential neighborhood $ U $ and a sequences of integers $ \{l_k\} $ such that the measures $ \{\nu^i_{l_k}\} $ converges in $*$-weak topology on $ H(U) $ for any $ i = 1, ..., 2deg(R) - 2.$
\endproclaim
\subheading{Acknowledgement}  I would like to
thank  to IMS at SUNY Stony Brook  FIM ETH at Zurich and IM UNAM at Cuernavaca for its hospitality during the preparation of this paper.

\heading{\bf Quadratic differentials for rational maps}\endheading
\par Let $ S_R $ be the Riemann surface associated with action of $ R $
on its Fatou set, then (see \cite{\bf S}) $ S_R $ is finite union $
\cup_iS_i $ of punctured torii  punctured spheres and foliated surfaces.
\par Let $ A(S_R) $ be space of quadratic holomorphic integrable
differentials on $ S_R $ and if $ S_R = \cup^N_iS_i, $ then $ A(S_R) =
A(S_1) \times ...\times A(S_N), $ where $ A(S_i) $ is the space of quadratic 
holomorphic integrable differentials on $ S_i. $
\subheading{Quadratic differentials for foliated surfaces} Due to Sullivan (\cite{\bf S})
a foliated surface $ S $ is either unit disk or round ring with marked points 
and is equipped with a  group
{ \bf G}$_D $ of rotations. This group  { \bf G}$_D $ is everywhere dense subgroup in 
the group of all rotations of $ S $ in the topology of uniform convergence on 
$ S. $ Hence for any $ z \in S $ the closure of the orbit  { \bf G}$_D(z) $ presents a
circle which is called {\it leaf of invariant foliation}. If the leaf $ l $ 
contains  a marked points $ x, $ then we call $ l $ as {\it critical leaf} and 
denote it by $ l_x. $ With exception of one case the boundary $ \partial S $ 
consists of critical leafs. This exception is the surface corresponding to the 
grand orbit of simply connected superattractive periodic component containing only one critical point. In this 
last case the surface $ S $ does not contain critical leafs. Note, that in this 
case we will assume that the modulus of $ S $ is not defined (see {\cite{\bf S}} for
 details).
 \par Any quadratic absolutely integrable 
 holomorphic differential $ \phi  $ has to be invariant under 
 action the group { \bf G}$_D $ for foliated surfaces.  Hence only $ \phi = 0 $ is the absolutely integrable holomorphic 
 function for $ S $ with undefined modulus and therefore we have in this case $ A(S) = \{0\}.$
\par
 After removing the critical leaves from $ S $ 
 we obtain the collection $ \cup S_i\cup D $ of the rings $ S_i $ and disk $ D $
 (in the case of Siegel disks). We call this decomposition as {\it critical 
 decomposition.} For this decomposition we have
 $$ 
 \phi_{S_i} =  h_i(z)\cdot dz^2, \text{ and } \phi_D = h_0\cdot dz^2
  $$
  where $ h_i $ is holomorphic absolutely integrable on $ S_i $ function and the 
  same for $ h_0 $ on $ D $. Easiest 
  calculations show that $ h_i(z) = \frac{c_i}{z^2} $ and $ h_0 = 0. $ 
  From the discussion above we conclude that for a ring with $ k $ critical 
  leaves (two from it present the boundary of $ S $) the dimension 
  $ dim(A(S)) = k - 1. $ 
  \par Now let $ S $ be ring with critical decomposition $ \cup_{i = 1}^k S_i$ and 
  $ \phi \in A(S) $ is a differential, then $ \Vert\phi\Vert = 4\pi\sum_i\vert 
  c_i\vert mod(S_i), $ where $ \phi = \sum_i\phi_{\vert S_i} = \sum_i 
  {\frac{c_i}{z^2}}_{\vert S_i} $ and $ mod(S_i) $ is modulus (or the extremal 
  length of the family of curves connecting the boundary component of $ S_i$) of the ring 
  $ S_i. $   
  \par We always assume here that the hyperbolic metric $ \lambda $ on the foliated 
  ring $ S $ is the collection of complete hyperbolic metrics $ \lambda_i $ on 
  components of critical decomposition of $ S. $ For example if $ \cup_i S_i $ is
  the critical decomposition of $ S, $ then the space $ H(S ) $ of harmonic 
  differentials on $ S $ consists of the elements 
  $$ \lambda^{-2}\overline{\phi} = \sum_i\frac{\overline{c_i}\lambda_i^{-2}}{\overline{z^2}}
  ,
   $$
  where $ \phi \in A(S). $
  \par The space of Teichmuller differentials $ td(S) $ consists of the elements
  $ \phi = 
  \sum_i {c_i\frac{z}{\overline{z}}\frac{d\overline{z}}{dz}}_{\vert S_i},
  $ where $ \cup_i S_i $ is the critical decomposition of $ S. $

\par
Denote by $ \Omega (R) $ the set 
$$ \overline{\Bbb C} \backslash \{\text{closure
of grand orbits of all critical points of }  R\}, 
$$
 then $ R $ acts on $
\Omega (R) $ as unbranched autocovering.
Let $ D \subset \Omega (R) $ be a 
periodic component, then $ D $ corresponds to either attractive or parabolic 
periodic domain.   
Let us fix the following
notations.
\roster
\item $ \Delta $ is unit disk,
\item $ j_D $ is universal covering $ \Delta \rightarrow D,$
\item let $ R^k: D\rightarrow D $ be the first return map for $ D, $ then 
let $ f $ be the  lifting of $ R^k $ onto 
$ \Delta $ by $ j_D. $ Note that $
f $ is a M\"obius map,
\item $ \Gamma $ is the group of the deck transformations of the
covering $ j_D $,
\item $ G = <f, \Gamma> $ is finitely generated Fuchsian group
uniformizing the surface $ S_D $ (i.e. $ \Delta / G \cong S_D $),
\item For a given Fuchsian group H let $ A(H) $ be the space of all
holomorphic functions on $ \Delta $ such that $ \phi(\gamma)(\gamma')^2 =
\phi $ for all $ \gamma \in H $ and for any $ \phi \in A(H) $ and
$$ \iint_\omega \vert\phi\vert < \infty $$
here $ \omega $ is a fixed fundamental domain for $ H. $ Further denote by $
B(H) $ the following space
$$\multline
 \{\phi \text{ holomorphic functions on } \Delta,
\phi(\gamma)(\gamma')^2 = \phi \text{ for all } \gamma \in H \text{ and
}\\
\sup_{z \in \Delta}\vert\lambda^{-2}\phi\vert < \infty, \text{ where }
\lambda \text{ is hyperbolic metric on } \Delta\},
\endmultline
$$
with norm $ \Vert\phi\Vert = \sup_{z \in
\Delta}\vert\lambda^{-2}\phi\vert $ the space $ B(H) $ presents a Banach
space. Now let $ S $ be a Riemann surface, then denote by $ B(S) $ the following
space
$$\multline
\{\phi \text{ quadratic holomorphic differentials on } S \text{ and }\\
\sup_{z \in S}\vert\lambda^{-2}\phi\vert < \infty, \text{ where }
\lambda \text{ is hyperbolic metric on } S\},
\endmultline
$$
with norm $ \Vert\phi\Vert = \sup_{z \in
S}\vert\lambda^{-2}\phi\vert $ the space $ B(S) $ presents a Banach
space.  If $ \Delta \big/_\Gamma \cong S, $ then it exists an
isometrical isomorphism $ \Phi $ from $ A(H) $ onto $ A(S).$
\endroster
\par Now let $ Y \subset \overline{\Bbb C} $ be an open subset. Then, as 
above, $ A(Y) $ denotes the space of holomorphic functions on $ Y $ 
absolutely integrable
over $ Y $ and $ B(Y) $ consists of holomorphic functions
$ \phi $ on $ Y $ with the following norm
$$ \Vert\phi\Vert = sup_{z \in Y}\vert\Lambda_Y^{-2}\phi\vert, $$
where $ \Lambda_Y $ is a  metric so that the its restriction over any component 
$ D \subset Y $ satisfies
$$ {\Lambda_Y}_{\vert D} = \lambda_D,
$$
where $ \lambda_D $ is Poincare metric on $ D.$ 
\proclaim{Lemma(Bers Duality Theorem)} 
\roster
\item The space $ B(H) $ ($B(S)$)
is isometrically isomorphic to the dual space $ A^\ast(H)
$ ($A^\ast(S)$) and this isomorphism is defined by the Peterson scalar
product
$$ \iint_\omega\lambda^{-2}\phi\overline{\psi}
(\iint_S\lambda^{-2}\phi\overline{\psi}), $$
where $ \omega $ is a fundamental domain for $ H $ and $ \phi \in
A(H)(A(S)), \psi \in B(H)(B(S)).$
\item If $ H $ is finitely generated group or the surface $ S $
 is compact surface with finite number punctures or $ S $ is foliated annuli. Then $
A(\Gamma) = B(\Gamma) $ and $ A(S)=B(S). $ Dimension of $ A(\Gamma) $ and $ (A(S) $ is finite and the Petersen scalar product becomes the inner scalar product.
\item Let $ Y \subset \overline{\Bbb C} $ be an open subset. Then
as above, the spaces $ A(Y) $ and $ B(Y) $ are dual by the Peterson 
scalar product
$$ \iint_Y(\Lambda_Y)^{-2}\phi\overline{\psi}.
$$
\endroster
\endproclaim
\demo{Proof} See \cite{\bf Kra}. \enddemo
\subheading{Poincare $\Theta-$operator for rational maps} We construct this
operator by the way which are similar to one in Kleinian group case. 
\par 1). Let $ D \in \Omega (R) $ is corresponding to an attractive periodic domain. 
\par Let $ \Theta_H(\phi) $ 
be theta series of Poincare for the
Fuchsian group $ H $, that is
$$ \Theta_H(\phi) = \sum_{\gamma \in H}\phi(\gamma)(\gamma')^2 $$
for $ \phi \in A(\Delta). $ This series defines the map from $ A(\Delta)
$ onto $ A(H) $ and the kernel is the space
$$ \text{ closure of the linear span}\{\phi - \phi(\gamma)(\gamma')^2,
\phi \in A(\Delta), \gamma \in H\} $$
\par In our case we have $ G = <f, \Gamma>. $ Let $ G = \cup_i \Gamma
g_i, $ then by $ \Theta_f $ we denote the relative $ \Theta-$ series
i.e.
$$ \Theta_f(\phi) = \sum_i\phi(g_i)(g_i')^2, $$
then it is clear that
$$ \Theta_G = \Theta_f \circ \Theta_\Gamma $$
Finally we define the $ \Theta-$series for our rational map $ R. $ Let $
\Psi $ be isometrical isomorphism from $ A(D) $ onto $ A(\Gamma), $
then we set.
$$ \Theta_D(\phi) = \Theta_f \circ \Psi(\phi), \text{ for } \phi \in
A(D). $$
So we have $ \Theta_D(A(D)) = A(G) \cong A(S_D).$
\par Define $ L(X) $ to be the {\it grand or full} orbit of the set 
$ X \subset \overline{\Bbb C}.$ Now we construct the map 
$ \Theta_{L(D)}:A(L(D))\rightarrow A(S_D) $ by the following way.
Let $ \phi \in A(L(D)) $ and $ X_i $ be components of $ L(D), $ then 
$ \phi = \sum_{i} \phi_{{\big\vert }X_i}. $
Let $ \Theta_{X_i}: A(X_i)\rightarrow A(D) $ be $ \Theta-$ operator corresponding to 
the unbranched covering $ R^{k_{X_i}} : X_i\rightarrow D, $ (where $ k_{X_i} $ is the minimal 
integer satisfying to the last property). Then we set
$$ 
 \Theta_{L(D)}(\phi) = 
 \Theta_D\left(\sum_{X \in L(D)}\Theta_X(\phi_{|X})\right)
$$
2) The case of parabolic domains $ D $ is similar to attractive one.\newline
3) Let $ D $ be a superattractive domain. This case corresponds to 
non-discrete
groups. Therefore we need an additional information respect to this foliated case.
Let us start with simple lemma about the Ruelle operator.  
\proclaim{Lemma 2} Let $ R $ be a rational map. Let $ Y \subset \overline{\Bbb C} $ be positive Lebesgue measure subset 
which is completely invariant under action of $ R, $
then the following is true.
\roster
\item $ R^* : L_1(Y) \rightarrow L_1(Y) $ is linear surjection with unit norm. The operator
$$ R_*(\phi)= \frac{\phi(R)(R')^2}{deg(R)}
$$ is an isometry "into" and $ R^*\circ R_* = I, $ where $ I $ is identity 
operator.
\item Beltrami operator 
$$
 B_R(\phi) = \phi(R)\frac{\overline{R'}}{R'}: L_\infty(Y) \rightarrow L_\infty(Y)
 $$
presents the dual operator to $ R^*. $  The operator $ B_R $ is an 
isometry "into".
\item If $ Y $ is an open set. Then $ R^*:A(Y) \rightarrow A(Y) $ is surjection 
of unit norm and $ R_* $ maps $ A(Y) $ into itself too. Let $ _*R: B(Y) 
\rightarrow B(Y) $ be the dual operator respect to Peterson scalar product. Then
$ _*R(\phi) = \phi(R)(R')^2 $ and $ _*R = \deg(R)R_*.$
\endroster
\endproclaim
\demo{Proof} All items are immediate consequences of the definition  of the operators.
\enddemo
\proclaim{Remark 3} If $ R: X \rightarrow Y $ is a branched covering for a rational 
map $ R $ and two domains $ X,Y \subset \overline{\Bbb C}.$ Then 
$ R^*: A(X) \rightarrow A(Y) $ is Poincare operator of the covering $ R. $
\endproclaim  

 The discussion above shows that the $ \Theta $ operator (in attractive and parabolic cases) is invariant respect to Ruelle operator, that is $ \Theta(R^*) = \Theta $ and hence $ ker(\Theta)\supset \overline{(I - R^*)(A(\Omega))}. $ 
\par Now let us continue the discussion on superattractive case. Let for simplicity $ D $ be an invariant superattractive domain. Our aim is to prove the following theorem. 
\proclaim{Theorem 4} Let $ D \subset F(R) $ be invariant superattractive domain and $ X = D \backslash Pc(R). $ Let $ S $ be the foliated surface associated with $ D .$ Then the quotient space $ A(X)_{{\big /}\overline{(I - R^*)(A(X)}} $ is isomorphic to the space $ A(S). $ 
\endproclaim
\demo{Proof} In further we need in the following basic facts about non-expansive operators.
\proclaim{Mean ergodicity lemma} Let $ T $ be non-expansive ($\Vert T\Vert =1$) linear endomorphism of a Banach space $ B $ and let  $ \phi \in B $ be any element.
\roster 
\item Assume that for {\rm Cesaro average} $A_N(T,\phi) = \frac{1}{N}\sum_{i =0}^{N - 1}T^i(\phi) $ there  exists
subsequence $ \{n_i\} $ such that $ A_{n_i}(T,\phi) $ weakly converges to
an element $ f \in B, $ then $ f $ is a fixed point for $ T $ and $ A_N(\phi)
$ converges to $ f $ strongly (i.e. by the norm). If $ f = 0 $ then $
\phi \in \overline{(I - T)(B)} $ and visa versa i.e. if $ \phi \in
\overline{(I - T)(B)}, $ then $ A_n(T,\phi) $ tends to zero with respect to  the norm.
\item The linear continuous operator $ T $ on a norm space $ B $ is
called mean ergodic if and only if the Cesaro average $ A_N(T,\phi)$  converges with respect to the  norm for any element $ \phi \in B.$ In this case $ B = Fix\times \overline{(I - T)(B)} $ and $ A_n(T,\phi) $ converges to projection $ P: B \rightarrow Fix, $ here $ Fix $ is the space of fixed elements for $T.$
\endroster
\endproclaim
\demo{Proof} See the book of Krengel (\cite{Kren}).
\par Now first consider the case of simply connectedness of $ D. $  Hence up to conformal changing of coordinates on $ D $ we can think 
that $ D = \Delta, X = \Delta^* = \Delta \backslash \{0\} $ and $ R(z) = z^2. $ 
\enddemo
In this case we claim
\proclaim{Claim} $ A(X) = \overline{(I - R^*)(A(X))}. $
\endproclaim
\demo{Proof of the claim} Let $ F $ be the space of finite linear combinations of the functions of the following kind $ \frac{1}{z - a}, $ where either $ a = 0 $ or $ \vert a\vert >1 $ or $ a $ is repelling periodic point for $ R(z). $ Then by Bers density theorem (see \cite{\bf Kra} or discussion below) we know that $ F $ presents everywhere dense subset of $ A(X).$
\par By using direct calculations we have.
$$ R^*\left(\frac{1}{z - a}\right) = \frac{a}{2z(z - a^2)} = \frac{1}{R'(a)}\left(\frac{1}{z} - \frac{1}{z - R(a)}\right)
$$
and hence by induction and the fact  $ R^*(\frac{1}{z}) = 0, $ we obtain
$$ (R^*)^n\left(\frac{1}{z - a}\right) = \frac{1}{(R^n)'(a)}\left( \frac{1}{z} - \frac{1}{z - R^n(a)}\right).
$$
Therefore for any $ \phi \in F $ we have $ (R^*)^n(\phi) \rightarrow 0, $ for $ n \rightarrow\infty $ exponentially fast. Hence Cesaro averages $ A_N(R^*, \phi) $ strongly tends to 0 for any $ \phi \in A(X). $ By Mean ergodicity lemma above we complete this claim.
\enddemo
\par Let us again consider unit disk $ \Delta $ and the map $ R(z) =z^2. $ Let $ b_1, ..., b_{k + 1} \in \Delta $ be points such that $ \vert b_{k + 1}\vert = \vert b_1\vert^2 < \vert b_k\vert < ... < \vert b_1\vert < 1. $   Let $ S $ be the ring $ \vert b_1\vert^2 \leq \vert z\vert \leq \vert b_1\vert $ with $ \{b_i\} $ as marked points and let $ G $ be the group of rotation of $ S. $ Then $ (S, G) $ is foliated surface.  Let $ X = \Delta \backslash\overline{\cup_iL(b_i)}, $ where like above $ L(b_i) $ means the grand orbit of $ b_i.$ We claim that.
\proclaim{Lemma 5} There exists continuous surjection $ P : A(X) \rightarrow A(S) $ such that $ ker(P) \supset \overline{(I - R^*)(A(X))}.$
\endproclaim
\demo{Proof} For simplicity assume that $ S $ has only two marked points $ x $ and $ R(x) $ on different components of $ \partial S. $ Let for any integer $ i $ the rings $ S_i, $ be component of $ L(S) $ such that $ S_0 = S $ and $ S_i = R(S_{i - 1}).$ Denote by $ W \in A(X) $ the subspace consisting of the elements $ \sum_i{\frac{c_i}{z^2}}_{{\big |}S_i}. $
\enddemo
 Then we claim.
\proclaim{claim} There exists projection $ Q : A(X) \rightarrow W $ such that
\roster
\item $ Q = lim P_n, $ where $ P_n = R_*^n\circ {R^*}^n $ and $ \Vert Q\Vert = 1.$
\item $ Q\circ R^* = R^*\circ Q $ and $ R_*\circ Q = Q\circ R_*, $
\item $ R^*, R_* : W \rightarrow W $ present surjective isometries.
\endroster
\endproclaim
\demo{Proof of the claim} Such as $ \Vert P_n \Vert \leq 1 $ and $ P_n(A(S_i)) \subset A(S_i) $ for any $n $ and $ i $ it is enough to show the convergence of $ P_n $ on $ A(S_i) $ for a fixed $i. $ Again by Bers density theorem the linear span of the function $ \frac{1}{z - a} $ with $ a \notin \overline{S_i} $ present everywhere dense subset of $ A(S_i).$
\par For any $ a \in {\Bbb C} \backslash \overline{S_i} $ we have
$$ P_n\left(\frac{1}{z - a}\right) = \frac{az^{2^n}}{z^2(z^{2^n} - a^{2^n})}
$$
and hence for big $ n > m $ obtain.
$$
\Vert (P_n - P_m)\left(\frac{1}{z - a }\right)\Vert = \iint_{S_i}\frac{\vert a\vert}{\vert z\vert^2}\left\vert\frac{1}{1 -
{\frac{\vert a\vert}{\vert z\vert}}^{2^n}} - \frac{1}{1 -{\frac{\vert a\vert}{\vert z\vert}}^{2^m}}\right\vert \vert dz\vert^2  = *
$$
 If  $ \vert a\vert = \lambda\cdot min_{z \in S_i}\vert z \vert $ with $ \vert \lambda\vert < 1, $ then
$$ * \leq 2 \lambda^{2^{m}}C\iint_{S_i}\frac{\vert a\vert}{\vert z\vert^2} \rightarrow 0 \text{ for } n > m \rightarrow \infty.
$$
If $ \vert a\vert = \lambda\cdot max_{z \in S_i}\vert z\vert $ with $ \vert \lambda\vert > 1, $ then
$$ * \leq \iint_{S_i}\frac{\vert a\vert}{\vert z\vert^2}\frac{\vert\left(\frac{a}{z}\right)^{-2m} - \left(\frac{a}{z}\right)^{-2n}\vert}{\vert \left(\frac{a}{z}\right)^{-2n} -1\vert\vert \left(\frac{a}{z}\right)^{-2m} - 1\vert} \rightarrow 0 \text{ when } n>m\rightarrow \infty.
$$
Let $ Q = \lim P_n $ in strong topology of $ L_1(S_i).$ Then we have\hfill\newline
1) if $ \vert a\vert > max_{z \in S_i}\vert z\vert, $ then $ Q\left(\frac{1}{z - a}\right) = 0 $ and\hfill\newline
2) if $ \vert a\vert < min_{z \in S_i}\vert z\vert, $ then $ Q\left(\frac{1}{z - a}\right) = \frac{a}{z^2}.$
\par 2) $ R^*\circ P_n = P_{n - 1} \circ R^* $ and $ P_n\circ R_* = R_*\circ P_{n - 1}, $ hence $ R^* $ and $ R_* $ commute with $ P. $
\par 3) By the construction  $ ker Q $ contains $ ker P_n $ for all $ n $ and hence $ R^*_{{\big\vert}W} $ is an isomorphism.
\par Now let us consider the map $ T: W\rightarrow A(S_0)=A(S) $ defined by the formula
$$ T(\phi) = \sum_{-\infty}^{\infty}(R^*)^n(\phi). 
$$
\enddemo
\proclaim{Lemma 6} $ T $ is continuous non-expansive operator with $ ker(T) \supset \overline{\left(I - R^*\right)(W)}. $
\endproclaim
\demo{Proof} Let $ \phi = \sum{\frac{c_i}{z^2}}_{{\big\vert}S_i} \in W, $
then $ T(\phi) = \frac{1}{z^2}\sum_{-\infty}^{\infty}c_i\cdot 2^i $ and
$ \Vert\phi\Vert = \sum\vert c_i\vert\iint_{S_i}\frac{1}{\vert z^2\vert} = 4\pi\sum\vert c_i\vert mod(S_i) = 4\pi\cdot mod(S_0)\sum 2^i\vert c_i\vert. $ Hence $ \Vert T(\phi)\Vert \leq \Vert\phi\Vert. $
\par The standard arguments imply that $ ker(T) \supset \overline{\left(I - R^*\right)(W)}. $ Lemma is complete.
\par By using  lemmas above we conclude that the map $ P = T\circ Q: A(X)\rightarrow A(S) $ is continuous surjection. 
\par Now let us return to theorem 4. Let $ h:U\subset D \rightarrow\Delta $ be conformal map with $ h(y) = 0 $ and $ h'(y) = 1, $ where $ y \in U $ is the fixed  superattractive point and $ h $ conjugates $ R $ with $ z\rightarrow z^2.$ Then the function $ \alpha = \log\vert h\vert $ can be harmonically extended onto $ L(D) $ to unique function, again denoted by $\alpha.$ Let $ c_i $ be critical points in $ L(D)\backslash y $ ordered by values $ \alpha(c_i). $ Then $ h $ can be conformally extended on region $ V = \{z, \alpha(z) < \alpha(c_1)\} $ and $ h(V) =\Delta_r=\{\vert z\vert < r < 1\}.$ 
\par The ring $ S_D = h(V)\backslash h(R(V)) $ with critical decomposition presents the foliated  surface associated with $ L(D). $
\par Let $ Y = L(D)\backslash \{\cup_iL(c_i)\cup L(y)\} $ and $ F = V\backslash \{\cup_iL(c_i)\cup L(y)\}, $ then as in attractive case construct the operator
$ \Theta_Y:A(Y)\rightarrow A(F).$   If $ h_*:A(F)\rightarrow A(X) $ is the injection generated by $h, $ where $ X $ is the set constructed by $ S_0 =S_D $ like in lemma 5. Then we set
$$ \Theta_{L(D)} = P\circ h_*\circ\Theta_Y: A(Y)\rightarrow A(S_D).$$
Theorem is proved.
\enddemo
\par Finally we set
$ \Theta(R): A(\Omega)\rightarrow A(S_R) $ by
$$ \Theta(R)(\phi) = \left(\Theta_{L(D_1)},..., \Theta_{L(D_k)}\right), 
$$
where $ D_i \subset F(R) $ are periodic components.
\enddemo
\subheading{Space $A(R)$}
\par Now again consider the space $ A(\Omega). $ Note that any function
of the kind
$$ \gamma_a(z) = \frac{a(a - 1)}{z(z - 1)(z - a)} \text{ for } a \in
\overline{\Bbb C} \backslash \Omega $$
belongs to $ A(\Omega).$  Let us introduce the subspace $ A(R) \subset
A(\Omega) $ as follows. 
Let $ S $ be the set 
$$ \{\cup_{i}\{L(c_i)\}\cup \{\{L(0, 1, \infty)\}
\backslash \{0, 1, \infty\}\}, $$
where $ c_i $ are critical points. Then we set
$$ A(R) = \text{linear span}\{\gamma_a(z), a \in S\} $$
This space $ A(R) $ is an linear space and we set on $ A(R) $ two different
topologies by the following norms $ \vert\cdot\vert_1 = \int_\Omega
\vert\cdot\vert $ and $ \vert\cdot\vert_2 = \int_{{\bold J}(R)}
\vert\cdot\vert $. Denote by $ A_i $ the spaces $ \{A(R),
\vert\cdot\vert_i\} $, respectively.
\proclaim{Remark 7} The space $ A(R) $ serves a kind of connection between spaces $ L_1(\Omega) $ and $ L_1(J(R)) $ and comparison of $\Vert\cdot\Vert_1 $ and $ \Vert\cdot\Vert_2 $ topologies is basis for our discussion below.
\endproclaim
\proclaim{Lemma 8} The operators $ R^*, R_* $ and $_*R $ are continuous endomorphisms of both spaces $ A_1 $ and $ A_2. $
\endproclaim
\demo{Proof} It is sufficient to show that for any $ \phi \in A(R) $ the functions $ R^*(\phi) $ and $R_*(\phi) $ belong to $ A(R) $ again.
\par Let $ \phi = \gamma_a. $ Then $ R^*(\phi) $ and $R_*(\phi) $ are holomorphic everywhere except finite numbers of points belonging to set $ S $ and hence are rational functions holomorphic on $\Omega. $ Besides both $ R^*(\phi) $ and $R_*(\phi) $ are integrable over $ \overline{\Bbb C} $ and hence belong to $ A(R). $ lemma is proved.
\enddemo

 Then we have the following
well-known result.
\proclaim{Lemma(Bers density theorem)} $ A_1 $ is everywhere dense in $
A(\Omega).$
\endproclaim
\demo{Proof} See for example book of I.Kra (\cite{\bf Kra}).
\enddemo
\proclaim{Lemma 9} Let $ L $ be a continuous functional on $ A_1 $
invariant under action of $ R^\ast $ (i.e. $L((R^\ast)(\phi)) =
L(\phi)$).Then $ L(\phi) = \iint\limits_\Omega
\lambda^{-2}{\overline{\psi}}\phi, $ where $ \lambda $ is hyperbolic
metric on $ \Omega $ and $ \psi \in B(\Omega).$
\endproclaim
\demo{Proof} Bers density and Bers duality theorems complete the proof of this lemma.
\enddemo
\heading{\bf  Bers Isomorphism} \endheading
\par Here we reproduce the Bers construction for the Beltrami
differentials and Eichler cohomology group with corrections (which
really often are evident) for the rational maps.
\par Consider the Beltrami action of $ R $ on the space $
L_\infty(\overline{\Bbb C}) $ i.e.
$$ B_R(\phi)(z) = \phi(R)(z)\frac{\overline{R'(z)}}{R'(z)}. $$
So the subspace $ Fix $ of fixed points for $ B_R $ in $ L_\infty(\overline{\Bbb
C}) $ is indeed the space of the invariant Beltrami differentials for $ R $
 unit ball of which describes all quasiconformal deformations of $
R.$
\par Now let $ {\bold J}_R $ be subspace of invariant Beltrami
differentials supported by { \bf J}($R$). Then it exists a continuous
map $ \Psi $ from $ A(S_R) \times {\bold J}_R $ into space $ Fix $ of fixed
points for $ B_R $ by the following way. Let $ \Theta^*:A^*(S_R) \rightarrow A^*(\Omega) $ be dual operator to $ \Theta-$operator. Then the image $ H(\Omega) = \Theta^*(A^*(S_R)) \subset B(\Omega) \subset Fix $ is called {\it space of harmonic differentials} and $ \dim(H(\Omega)) = \dim(A^*(S_R)) = \dim(A(S_R)).$ Let $ \alpha: A(S_R)\rightarrow A^*(S_R) $ be isomorphism defined by Petersen scalar product. Then we can define
$$ \Psi: A(S_R)\times {\bold J}_R \rightarrow Fix \text{ by } \Psi(\phi, \mu) = \left(\Theta^*\circ\alpha(\phi), \mu\right).
$$
\par Now normalize $ R $ so that 0, 1, $\infty $ are fixed points for
$R.$ Let $ C $ be component of the subset of rational maps in $ {\Bbb C}P^{2d + 1} $
fixing the points 0, 1 and $ \infty $ containing $ R. $ By $ H^1(R) $ we
denote the tangent space to $ C $ at the point $ R.$ Then $ H^1(R) $ may
be presented as follows. Let $ R(z) = z\frac{P_0}{Q_0}, $ then
$$ H^1(R) = \{z\frac{PQ_0 - QP_0}{Q^2_0}, \text{ where } Q(1) =
P(1), deg(Q) \leq deg(R), deg(P) \leq \deg(R) - 1\}, $$
where $ P, Q $ are polynomials and $ dim(H^1(R)) = 2d - 2.$
\proclaim{Remark 10} We use the notation $ H^1(R) $ by the following
reasons
\roster
\item The Weyl cohomology' construction for the action of
$ R $ (by the formula $ {\tilde R}(f) = \frac{f(R)}{R'}$) on the space  of
all rational functions gives the space $ H $ which is isomorphic to tangent
space to $ {\Bbb C}P^\infty $ at $ R $ up to normalization. More
precisely $ H $ is equivalent to direct limit
$$ H^1(R)\overset{j_1}\to\longrightarrow H^1(R^2)\overset{j_2}\to \longrightarrow H^1(R^3) ..., $$
where $ j_i $ are equivalent to the action $ \tilde R. $
\item This construction for Kleinian group is called Eichler
cohomologies.
\endroster
\endproclaim
\par Now follow Bers (see for example \cite{\bf Kra}) we introduce
 the Bers map $ \beta $ from $ {\bold J}_R \times A(S_R)
$ into $ H^1(R). $
\par Let $ \mu \in L_\infty(\Bbb C),$ then the function
$$ F_\mu(z) = z(z - 1)\iint\limits_{\Bbb C}\frac{\mu d\xi
d{\overline{\xi}}}{\xi(\xi - 1)(\xi - z)},$$
is continuous on $ \Bbb C $ and $ \vert F(z)\vert = O\vert
z^2\vert $ for $ z \rightarrow \infty.$ 
$$ \frac{\partial F_\mu}{\partial {\overline z}} = \mu. $$
in sense of distribution. $ F_\mu $ is called {\it potential } for $\mu.$
\par Let us define the Bers map $ \beta(t) $ for $ t = (\mu, \phi) \in
 A(S_R)\times {\bold J}_R $ by the formula
$$ \beta(t = (\phi, \mu))(z) = F_{\Psi(\phi, \mu)}(R(z)) -
R'(z)F_{\Psi(\phi, \mu)}(z). $$
\proclaim{Theorem 11}
\roster
\item $ \beta $ is injective antilinear map,
\item $ \beta({\bold J}_R \times A(S_R)) \subset H^1(R), $
\item if $ R $ is structurally stable, then $ \beta $ is an isomorphism.
\endroster
\endproclaim
\demo{Proof} Fix $ t \in {\bold J}_R \times A(S_R), $ then on $ {\Bbb C}
\backslash \{\text{poles of R}\} $ the derivative $
\frac{\partial\beta(t)(z)}{\partial{\overline z}} $ (in sense of distribution)
is zero, hence by Weyl lemma $ \beta(t)(z) $ is
holomorphic on $ {\Bbb C} \backslash \{\text{poles of R}\}. $ Further,
poles of $ R $ are at most than poles for $ F $ and we conclude that $
\beta(t)(z) $ is a rational function with zeros at 0 and 1 and simple
pole at $ \infty. $
\par Let k be the norm $ \Vert\Psi(t)\Vert_{L_\infty(\Bbb C)}. $
Consider the disk of Beltrami differentials $ \mu_x(z) =
x\Psi(t)(z), $ for $ \vert x\vert < k. $ let $ f_x $ quasiconformal maps
normalized by $ f_x(0, 1, \infty) = 0, 1, \infty, $ respectively. Then
for small $ x $ we have
$$ f_x(z) = z - z(z - 1)\iint\limits_{\Bbb C}\frac{x\Psi(t)(\xi)d\xi
d{\overline{\xi}}}{\xi(\xi - 1)(\xi - z)} + O(\Vert
x\Psi(t)\Vert^2_{L_\infty({\Bbb C})}),$$
and hence
$$ {\frac{\partial f_x(z)}{\partial x}}\vert_{x=0} = -F_{\Psi(t)}(z).$$
But $ R_x(z) = f_x \circ R \circ f^{-1}_x(z) $ are rational maps and so
by differentiation respect to $ x $ of equality above one have
$$ \frac{\partial f_x}{\partial x}\vert_{x = 0}(R(z)) -
R'(z)\frac{\partial f_x}{\partial x}\vert_{x=0}(z) = \frac{\partial
R_x}{\partial x}\vert_{x=0}(z) \in H^1(R).$$
Hence we have
$$ F_{\Psi(t)}(R(z)) - R'(z)F_{\Psi(t)} = - \frac{\partial
R_x}{\partial x}\vert_{x=0}(z) \in H^1(R).$$
\par Now for finishing lemma it is enough to assume that $ R $ does not
have conformal centralizer and use the Sullivan result (see \cite{\bf
S}) which particularly says that $ dim(T(R)) = dim(T_{R}({\bold c}(R)), $ where $
T_{R}({\bold qc}(R)) $ is tangent space to $ {\bold qc}(R) $ in initial
point and $ T(R) $ is Teichmuller space of $ R. $ Besides we
have $ dim(H^1(R)) = dim(T_{{\bold qc}(R)}(R)) $. So theorem is proved.
\enddemo
\heading{\bf Beltrami Differentials on Julia set} \endheading
\par Here we consider in details the space $ {\bold J}_R. $ Each element
$ \mu \in {\bold J}_R $ defines an $ {\underline invariant} $ respect to
the Ruelle operator functional $ L_{\mu} $ on the space $ A(R) $ which is continuous
in topology of the space $ A_2 $ (we recall that $ A_i = (A(R),
\vert\cdot\vert_i)).$ Continuity of $ L_\mu $ in the topology of the
space $ A_1 $ is crucial in the question of non-triviality of $ \mu.$
Indeed we have the following lemma.
\proclaim{Lemma 12} Let $ \mu \in {\bold J}_R, $ then $ \mu = 0 $ if
and only if $ L_\mu $ is continuous functional on $ A_1.$
\endproclaim
\demo{Proof} Let $ L_\mu $ is continuous on $ A_1 $. Then $ L_\mu $ is
continuous on $ A(\Omega) $ by the density theorem. Then by  lemma
9 the functional $ L_\mu $ may be presented by the expression $
L_\mu(\alpha) = \iint\limits_{\Bbb C}\alpha\lambda^{-2}\overline{\psi},
$ for some $ \psi \in A(S_R) $ and hence $ F_\mu(a) = L_\mu(\gamma_a) =
\iint\limits_{\Bbb C }\gamma_a\lambda^{-2}\overline{\psi} =
F_{\lambda^{-2}\overline{\psi}}(a). $ In other words we have
$ \beta(t = (\mu, 0))(z) = \beta(0, \psi)(z) $ for $ z \in \overline{S}, $
but $ \beta(t)(z) $ is rational map for any $ t \in A(S_{R}) \times
{\bold J}_R $ and $ \overline{S} $ is closed infinite set, hence we
conclude
$$ \beta(t = (\mu, 0)) = \beta(0, \psi) $$
that contradicts with injectivity of $ \beta.$ Lemma is proved.
\enddemo
\par Now we begin to consider the relations between continuity of $
L_\mu $ for $ \mu \in {\bold J}_R $ and properties of the Ruelle
operator $ R^\ast : A_2 \rightarrow A_2. $ Recall that operator $ R^\ast $
acts as linear autosurjection of $ L_1({\bold J}(R)) $ with unit norm.
\proclaim{Proposition 12} Let $ R $ be a rational map with simple critical points.
Assume that Lebesgue measure of $ Pc(R)$ is zero. Then $ {\bold J}_R =
\emptyset $ if and only if the Ruelle operator $ R^\ast : A_2
\rightarrow A_2 $ is mean ergodic.
\endproclaim
\demo{Proof} If $ {\bold J}_R = \emptyset, $ then subspace $ (I -
R^\ast)(A_2) $ is everywhere dense in $ A_2 $ and by the item (2)
of Mean ergodicity lemma we complete proof.
\par Now suppose that $ R^\ast $ is mean ergodic on $ A_2. $ Let $ \mu
\neq 0 \in {\bold J}_R, $ then by  lemma 12 $ L_\mu $ is not
continuous functional on $ A_1.$ But $ (I - R^\ast)(A_1) \subset
ker(L_\mu) $.
\par Now we claim that
$$ dim(A_1 /((I -R^\ast)(A_1))) < \infty $$
\demo{Proof of the claim} We start with the following lemma. Let us consider an element $ \tau_a  = \frac{1}{z - a} $ with $ a
\in {\Bbb C}. $ Now we show that
\proclaim{Lemma 14} For any integer $ n $ we have
$$
 {R^\ast}^n(\tau_a)(z) = \frac{1}{(R^n)'(a)(z - R^n(a))} - \sum_i\frac{b^n_i}{(a - c^n_i)(z - R(c^n_i))},
$$
where $ b^n_i $ are coefficients from the following decomposition $ \frac{1}{(R^n)'(z)} = \sum_i\frac{b^n_i}{z - c^n_i} + p_n(z) $ the points $ c^n_i $ are critical points of $ R^n $ and $ p_n(z) $ is polynomial.
\endproclaim
\demo{Proof of lemma} Let $\phi $ be a differentiable function with compact support, then for any fixed $ n $ one obtain
$$ 
\multline
\iint\phi_{\overline z}{R^*}^n(\tau_a)(z)=\iint\phi_{\overline z}(R^n)\frac{\overline{(R^n)'}}{(R^n)'}\tau_a=\iint\frac{(\phi\circ R^n)_{\overline z}}{(R^n)'}\tau_a= 
\sum_ib^n_i\iint\frac{(\phi\circ R)_{\overline z}}{(z - c_i)(z - a)}
\\+\iint\frac{(\phi\circ R)_{\overline z}\cdot p_n(z)}{z - a} = \sum_i\frac{b^n_i}{(a - c^n_i)}\left(\iint\frac{(\phi\circ R^n)_{\overline z}}{z - a} - \iint\frac{(\phi\circ R^n)_{\overline z}}{z - c^n_i}\right) +\\
+ \iint\frac{(p_n(z)\cdot\phi\circ R)_{\overline z}}{z - a} 
= \frac{\phi(R^n(a))}{(R^n)'(a)} - \sum\frac{b^n_i}{a - c^n_i}\phi(R^n(c^n_i)) =\\
= \iint\phi_{\overline z}\left(\frac{1}{(R^n)'(a)(z - R^n(a))}-  \sum_i\frac{b^n_i}{(a - c^n_i)(z - R^n(c^n_i))}\right).
\endmultline
$$
Hence by Weyl lemma we have that function
$$ h_n(z) = {R^\ast}^n(\tau_a)(z) - \frac{1}{(R^n)'(a)(z - R^n(a))} - \sum_i\frac{b^n_i}{(a - c^n_i)(z - R^n(c^n_i))}
$$
is integrable and holomorphic on $ \overline{\Bbb C} $ and hence $ h(z) = 0.$
\enddemo
\proclaim{Remark 15} Assuming that $ z=\infty $ is fixed point for $ R, $ we easily calculate that $ p_n(z) = \lambda^n, $ where $\lambda $ is multiplier of $\infty. $ For example, if $ F(R) $  has an attractive component, we always can think (by qc-surgery) that points $ z = \infty $ is a superattractive point.
\endproclaim
\par Now by using lemma 14 and the fact $ \gamma_a = (a - 1)\tau_o - a\tau_i + \tau_a $ we can write
$$ R^\ast(\gamma_a) = \sum_i \omega_i\gamma_{R(c_i}(z) + \omega\gamma_{R(a)}(z)
\tag{*}.$$
So we have
$$ R^\ast(\gamma_a) = \gamma_a - (\gamma_a - R^\ast(\gamma_a)) = \sum_i
\omega_i\gamma_{R(c_i)}(z) + \omega\gamma_{R(a)}(z). $$
If $ a \in R^{-1}(0, 1, \infty), $then again the direct calculation
shows that
$$ R^\ast(\gamma_a)(z) = \sum_i\alpha_i\gamma_{R(c_i)}(z).$$
We conclude that for any element $ \phi \in $ linear span$\{\gamma_a; a \in
S \backslash \{$the forward orbits of critical values$\}\} $ it exists $
n $ such that $ {R^\ast}^n(\phi) \in $linear span$\{\gamma_a; a\in
\{$forward orbits of critical points$\}\}. $
\par Now let $ a = R^k(b), $ where $ b $ is a critical value, then
$ \gamma_b \thicksim {R^\ast}^k(\gamma_b)
= \sum_i\alpha_i\gamma_{b_i} + \alpha\gamma_a $ that means
$ \alpha\gamma_a \thicksim \gamma_b -
\sum_i\alpha_i\gamma_{b_i}, $ (here $ \phi \thicksim \psi $ iff $ \phi
- \psi \in (I - R^\ast)(A_1)$).
\par As result we obtain that the space $ A_1 /((I -R^\ast)(A_1)) $ is
isomorphic to a subspace in $ X = $linear span$\{\gamma_{R(c_i)}\}.$
Claim is proved.
\enddemo
\proclaim{Remark 16} Note that finiteness of $ dim(A_1 /((I -R^\ast)(A_1))) $ is purely
algebraic fact does not depending on topology.
\endproclaim
\par By the claim we conclude if $ ker(L_\mu) $ contains closure of $ (I
- R^\ast)(A_1) $ in $ A_1, $ then $ L_\mu $ is continuous.
\par We claim that $ ker(L_\mu) $ {\it contains the space} $ \overline{(I -
R^\ast)(A_1)}.$
\demo{Proof of the claim} Otherwise it exists an element $ \phi \in \overline{(I -
R^\ast)(A_1)} $ such that $ L_\mu(\phi) \neq 0. $ By Mean ergodicity lemma we have that Cesaro averages $ A_N(\phi) $
tends to zero in $ A_1 $ by the norm. Further by
invariance of $ \mu $ we have
$$ L_\mu(A_N(\phi)) = \iint\limits_{{\bold J}(R)}\mu A_N(\phi) =
L_\mu(\phi). $$
Now using the mean ergodicity of $ R^\ast $ we obtain convergence of $
A_N(\phi) $ to an element $ f \in L_1({\bold J}(R)) $ by the norm and
hence $ L_\mu(\phi) = L_\mu(A_N(\phi)) = \iint\limits_{{\bold
J}(R)}\mu f.$ Besides the functions $ \{A_N(\phi)\} $ forms the normal
family of holomorphic functions on
$ Y = \overline{\Bbb C} \backslash \bigl\{ \{$closure of forward orbit of a$\}
\cup\{Pc(R)\}\bigr\}. $ Such that $ \Vert A_N(\gamma_a)\Vert_{A_1}
\rightarrow 0 $ we obtain that $ A_N(\phi) $ converges to zero
uniformly on compacts in $ Y $. Hence $ f = 0 $ on $ Y $ and such as $
m({\Bbb C} \backslash Y) = 0 $ one has $ L_\mu(\phi) = 0.$
Contradiction. Claim is proved.
\enddemo
\par So $ L_\mu $ is continuous functional on $ A_1 $ and by lemma
12 $ \mu = 0.$ Proposition is proved.
\proclaim{Definition} We call a rational map {\rm mean ergodic } if and only if the Ruelle operator  $ R^\ast : A_2
\rightarrow A_2 $ is mean ergodic.
\endproclaim
\enddemo

Now we show that topologies $ \Vert\cdot\Vert_1 $ and $ \Vert\cdot\Vert_2 $ are "mutually disjoint."
Denote by $ X_i $ the closure of the space $ \left(I - R^*\right)\left(A(R)\right) $ in spaces $ A_1 $ and $A_2 .$
\proclaim{Proposition 17} Let $ R $ be a rational map and $ dim(A(S_R))\geq 1. $ Then the following conditions are equivalent.
\roster
\item the map $ i=id : A_1 \rightarrow  A_2 $ maps weakly convergent sequences onto weakly convergent sequences.
\item $ i(X_1) \supset X_2, $ 
\item the Lebesgue measure of Julia set is zero.
\endroster
\endproclaim
\demo{Proof} Condition (3) trivially implies the conditions (1) and (2). 
\par Assume condition (1) is hold. then the dual map $ i^*: A^*_2 \rightarrow  A^*_1
 $ is continuous in $\ast-$weak topologies on $ A^*_1 $ and $ A^*_2. $ Hence for any
  $ \mu \in A^*_2 \subset L_\infty(J) $ there exists an element $ \nu \in A^*_1
   \subset L_\infty(F) $ such that $ \nu = i^*(\mu) $ and
$$ \iint_J\mu\gamma = \iint_F\nu\gamma.$$
Then for any $ \gamma \in A(R) $ we have $\iint_{\Bbb C}\gamma (\mu - i^*(\mu)) = 0. $ Let $ F_\mu(z) $ and $ F_\nu(z) $ be potentials. Then $ F_{{\big\vert}J(R)} =  \left(F_\mu(z)- F_\nu(z)\right)_{{\big\vert}J(R)} = 0 $ and if $ m(J(R) > 0 $ we have $ F_{\overline{z}} = 0 $ almost everywhere on $ J(R), $ where $ F_{\overline{z}} $ in sense of distributions. Hence we deduce:
$$ \mu - i^*(\mu) = 0 $$
almost everywhere on $ J(R). $ Since $F(R)\cap J(R) =\emptyset $ we have $ \mu = 0  $ almost everywhere and  we conclude that $ A^*_2=\{0\}. $ Hence $A_2 =\{0\} $ and we obtain $ m\left(J(R)\right) = 0.$ 
\par Now assume (2). Then conditions imply that any invariant continuous functional on $ A_1 $ generate an invariant line field on Julia set that contradicts to injectivity of the Bers map. Now using the fact that we can always think that $ F(R) $ contains an attractive domain we conclude $ m\left(J(R)\right) =0. $ Proposition is proved.
\enddemo
\proclaim{Proposition 18} Assume that $ F(R) \neq \emptyset $ and $ m\left(J(R)\right) > 0 $ for the given rational map $ R.$ Then 
there exists no invariant line fields on Julia set if and only if $ i^{-1}(X_2) \supset X_1.$
\endproclaim
\demo{Proof} It is easy. If there exists no invariant line field, then $ X_2 = A_2.$ Now assume $ i^{-1}(X_2) \supset X_1,$ then existing of invariant line field contradicts to injectivity of the Bers map.
\enddemo
\par We finish this chapter with simple application of discussion above to polynomials of degree two.
\proclaim{Theorem B} Let $ R(z) = z^2 +c $ and $ S_L = \sum^{L}_{j = 0}\frac{1}{(R^{j})'(c)} $ Assume there exists a subsequence $\{n_i\} $ of integers such that
\roster
\item $ \lim_{i \rightarrow \infty} \vert (R^{n_i})'(c)\vert =\infty $ and
 $\overline{\lim}_{i \rightarrow \infty} \vert S_{n_i}\vert > 0 $ or
\item $ \vert (R^{n_i})'(c)\vert \sim C=Const $ for  $ i\to\infty $ and 
$\overline{\lim}_{i \rightarrow \infty} \vert S_{n_i}\vert = \infty $ 
\endroster
Then there exists no invariant line field on Julia set.
\endproclaim
\demo{Proof}Our aim is to show that under conditions  $ X_2 = A_2.$ We know that $ \dim\{{A_2}_{\big/{X_2}}\} = 1.$ Hence any element $ \gamma \in A_2 $ has the following expression
$$ \gamma(z) = A\frac{1}{z - c} + \gamma_1, \text{ where } \gamma_1 \in X_2. $$
Let us show that element $ \gamma_c =  \frac{1}{z - c} $ belongs to $ X_2. $
Let us define the sequences
$$ \phi_0 = \gamma_c, \phi_i = \gamma_c + \sum_{j = 1}^{n_i - 1}\frac{1}{(R^j)'(c)}\frac{1}{z - R^j(c)}, $$
then 
$$ \phi_i - R^*({\phi_i}) = \frac{1}{z -c}\left(\sum_{j = 0}^{n_i}\frac{1}{(R^j)'(c)}\right) - \frac{1}{(R^{n_i})'(c)}\left(\frac{1}{z - R^{n_i}(c)}\right). 
$$
 Now assume that $ \gamma_c \notin X_2, $ then there exists a linear functional $ f $ on $ A_2 $ such that $ f(\gamma_c) \neq 0 $ and $ f_{|X_2} = 0. $ Hence $ f(R^*(\gamma)) = f(\gamma) $ and we calculate
$$ 0 = f(\phi_i - R^*(\phi_i)) = f(\gamma_c)\sum_{j = 0}^{n_i}\frac{1}{(R^j)'(c)} - \frac{1}{(R^{n_i})'(c)}f\left(\frac{1}{z - R^{n_i}(c)}\right). $$
By using condition we conclude that $ f(\gamma_c) = 0.$ Proposition 18 completes theorem B.
\enddemo
\heading{\bf Convergent rational maps}
\endheading
 \par We start with accumulation some facts (see books of I.Kra
 "Automorphic forms and Kleinian Group" I. N. Vekua "Generalized
 analytic function.")
\proclaim{Facts} Denote by $ F_{\mu}(a) $ the following integral $ \iint_{\Bbb
 C}\mu(z)\tau_a(z)dzd{\overline{z}} $ where $\tau_a(z) = \frac{1}{z -
 a}$ for $ a \in {\Bbb C}$ and $ \mu \in L_{\infty}(J(R)). $
 Then
 \roster
 \item $ F_{\mu}(a) $ is continuous function on $ {\Bbb C} $ and
 $\frac{\partial F_{\mu}(a)}{\partial {\overline{z}}} = \mu $ in sense
 of distributions.
 \item $\mid F_{\mu}(a)\mid = O(\mid z\mid^{-1}) $ for big $ z.$
$\Vert F_{\mu}(a)\Vert_{\infty} \leq \Vert\mu\Vert_{\infty}M, $
where $ M $ does not depends on $ \mu $ and $ a \in {\Bbb C} $.
\item $\vert F_{\mu}(a_1) - F_{\mu}(a_2)\vert \leq
\Vert\mu\Vert_{\infty}C\vert a_1 - a_2\vert\vert\ln\vert a_1 -
a_2\vert\vert $, where $ C $ does not depends on $ \mu $ and $ a $.
\endroster
\endproclaim
\par  Denote by $ B $ the operator $ \mu \to F_{\mu}(a) : L_\infty(J(R))
\to C(\Bbb C) $ and by $ X $ the image $ B(L_\infty(J(R)) $.
Now by $ W $ denote the space $ X $ with the following topology:
$$ \phi_n \to 0 \text{ means } \Vert\phi_n\Vert_{\infty} \to 0
\text{ and } \frac{\partial\phi_n}{\partial\overline{z}} \to 0 \text{ in
$*$-weak topology of } L_\infty(J(R)).
$$
\proclaim{Lemma 19}
\roster
\item $ W $ is complete locally convex vector topological space.
\item B is the compact operator mapping $L_\infty(J(R))$ onto $ W $,
that is $ B $ maps bounded sets onto precompact sets.Here precompactness
means any sequence contains "Cauchy" subsequence.
\item Any bounded set $ U \subset W $ is precompact.
\endroster
\endproclaim
\demo{Proof} The first is evident.
\par 2). Let $ U \subset L_{\infty}(J(R)) $ is bounded, then $ U $ is
precompact in $*$-weak topology of $ L_{\infty}(J(R)) $. Further from
item 2 of Facts we have that $ B(U) $ forms uniformly bounded and
equicontinuous family of continuous functions that means $ B(U) $ is
precompact in topology of uniform convergence.
\par 3). What that means boundedness in $ W $? Particularly the set
$$ V = \{\frac{\partial\phi}{\partial\overline{z}} \text{ in sense of
distributions, for }\phi \in U \}
$$
forms bounded set in $*$-weak topology of $L_{\infty}(J(R)) $ hence $
V $ is bounded in the norm topology of $L_{\infty}(J(R)) $. We complete
 lemma by using the item 2 and the fact $ \phi =
B(\phi_{\overline{z}}).$
\enddemo
\par Let us define the operator $ T $ on $ X $ as follows
$$
T(F_\mu(a)) = F_{B_R(\mu)}(a) = \iint_{\Bbb C}B_R(\mu)\tau_a =\iint_{\Bbb C}\mu R^*(\tau_a(z)),
$$
where $ B_R $ is the Beltrami
operator. Easily show that
$$
T(\phi) = \frac{\phi(R(a))}{R'(a)} - \sum\frac{b_i\phi(R(c_i))}{a -
c_i},
$$
where $ \sum\frac{b_i}{a - c_i} = \frac{1}{R'(a)}.$
For example for $ R(z) = z^2 +c $ we have $ T(\phi)(a) = \frac{\phi(R(a)) - \phi(c)}{R'(a)}.$
\proclaim{Remark 20} From definition we see that
$$ \{T^n(\phi), n = 0, 1, ...\}
$$
forms bounded in $ W $ set
\endproclaim
\proclaim{Lemma 21} $ T $ is continuous endomorphism of $ W $.
\endproclaim
\demo{Proof} Let $ F_{\mu_i} \to 0 $ in $ W $, then $ \Vert\mu_i\Vert
\leq C < \infty $ and hence $\{T(F_{\mu_i})\}$ forms precompact in $
W $  family. Let $ \psi_0 $ be a limit point of this set, then
$$ \psi_0(a) = \lim_j T(F_{\mu_{i_j}}) = \iint_{\Bbb
C}\mu_{i_j}R^*(\tau_a) \to 0 (*-\text{weak topology}).
$$
So $ \psi_0 = 0.$ \enddemo
\par Now we start with weak conditions implying the mean ergodicity of given rational map $ R. $
\proclaim{Definition} We will say that a rational map $ K(z) $ is {\rm
strongly convergent} if the space of
 $ qc_J(K(z)) $ contains a map $ R $ for which
there exists a point $ d $ with 
$$
card(\cup_0^\infty R^n(d)) >
2d - 1 \text{ and
}
 s_n(d) < \infty  \text{ for all } n = 1, ...,
$$
 where
$$
 s_n(z) = \sum\frac{\vert b^n_i\vert}{\vert z - c^n_i\vert} \tag{$*$}$$

and here $ \sum\frac{b^n_i}{z - c^n_i} + p_n(z) = \frac{1}{(R^n)'(z)}.$
\endproclaim

\par Now and below we assume that $ R(z) $ is the map from condition
$\ast $ with point $ z=\infty $ as superattractive point.
\par First note that
$$ s_n(R^m(a)) \leq s_{n + m}(a)\vert (R^m)'(a)\vert.$$
Indeed
$$ \frac{s_n(R^m(a))}{\vert (R^m)'(a)\vert} =
\frac{1}{\vert (R^m)'(a)\vert}\sum\frac{\vert b_i\vert}{\vert R^m(a) -
c_i\vert} = \sum\frac{\vert b_l\vert}{\vert a - c_l\vert}\sum\frac{\vert
Q(a)b_i\vert}{\vert P(a) - Q(a)c_i\vert} = ** $$
here $ R(a) = \frac{P(a)}{Q(a)} $
$$ ** \leq \sum\frac{\vert \gamma_k\vert}{\vert a - c_k\vert}  = s_{n +
m}(a).$$

\proclaim{Lemma 22} Assume $ s_n(a) \leq C < \infty $ for all $ n $ for
the given rational map $ R. $ Then Cesaro averages $ A_N(\tau_x) $
converges with the $ L_1(J)-$norm for any $ x \in \cup_{l=0}^\infty R^l(a).$
\endproclaim
\demo{Proof} In notations of above we have
$$\split
T(F_\mu)(y)
&= \iint_J\frac{B(\mu)}{z - y} dz\land d{\overline{z}} =
\iint_J\mu R^\ast(\tau_y)dz\land d{\overline{z}}\\
&= \frac{F_\mu(R)(y)}{R'(y)} - \sum\frac{b_i F_\mu(R(c_i))}{y - c_i} =
\sum\frac{b_i(F_\mu(R(y)) - F_\mu(R(c_i)))}{y - c_i}.\endsplit
$$
Now consider the sequence of functionals $ l_i(F) = (A_i(T)(F))(a) $ on
$ W. $ Under assumption we have
$$ \vert l_i(F)\vert \leq 2\frac{1}{i}\sum_{j=0}^{i -1}s_j(a)\sup_{w \in {\Bbb C}}\vert F(w)\vert $$
and so the family functionals $\{l_i\} $ can be continued on the space
$ C(\overline{\Bbb C}) $ of continuous functions on $ \overline{\Bbb C}
$ to family of uniformly bounded functionals. Therefore we can choose a
subsequence $ l_{i_j} $ converging pointwise to some continuous
functional $ l_0 $ (Note that $ l_0 $ is the fixed point for the dual
operator $ T^* $ acting on dual $ W^*$). Besides that means sequence
$ A_{i_j}(R^*)(\tau_a) $ weakly converges in $ L_1(J) $ and hence by the
standard ergodic arguments whole sequence $ A_N(\tau_a) $ converges by
the norm to a fixed for $ R^* $ element.
\par By notation above we know that
$$ s_n(x) \leq C_x, \text{ for any } x \in \cup R^n(a).$$
So by repeating of the arguments above we complete proof.
\enddemo
\par Now we prove main theorem of this chapter.
\proclaim{Theorem 23} Let condition $ \ast $ holds for the rational map $
R $, then  $ R $ is mean ergodic. 
\endproclaim
\demo{Proof} Proof of theorem we divide onto two steps.
\par The first step consists of proving  theorem under additional
assumption. Namely, Let $ d_i, i = 1,..., k $ be all critical values of $ R $, now
form $ k $ families of functionals on $ W $ like in lemma above $
l^i_n(F) = A_n(T)(F)(d_i) $. Now assumption is:
{\it It exists a subsequence $ \{n_m\} \subset \{n\} $ so that for all $
i \leq k $ subsequences $ l^i_{n_m} $ are convergent pointwise on $ W
$.}
\par We know that $ A $ is the closure of the linear span of the family
of functions $ \{\frac{a(a - 1)}{z(z - 1)(z - a)}, a \in S\}.$ So it is
enough to show convergence $ A_N(R^*)(\frac{a(a - 1)}{z(z - 1)(z - a)})
$ for any fixed $ a \in S.$
\par Now let $ x \in J(R) $ be any periodic repulsive point, then we
{\bf claim} that the sequence of functionals $ L_{n_m}(F) = A_{n_m}(T)(F)(x) $
converges pointwise on $ W $.
\par Proof of the claim. Without loss of generality we assume that $ x $
is fixed point for $ R. $ Let us calculate
$$ T(F)(a) = \frac{F(R(a))}{R'(a)} - \sum \frac{b_i F(R(c_i))}{a -
c_i},$$
$$ T^2(F)(a) = \frac{F(R^2(a))}{(R^2)'(a)} - \frac{1}{R'(a)}\sum_i
\frac{b_i F(R(c_i))}{R(a) - c_i} - \sum_i\frac{b_i T(F)(R(c_i))}{a - c_i},
$$
and
 
$$ \multline
T^n(F)(a) = \frac{F(R^n(a))}{(R^n)'(a)} - \frac{1}{(R^{n-1})'(a)}\sum_i
\frac{b_i F(R(c_i))}{R^{n-1}(a) - c_i} -\\
-\frac{1}{(R^{n-2})'(a)}\sum_i\frac{b_i T(F)(R(c_i))}{R^{n-2}(a) - c_i} - ...
- \sum_i\frac{b_iT^{n-1}(F)(R(c_i))}{a - c_i}.
\endmultline
$$
So for $ a = x $ and $ \lambda = R'(x) $ we conclude
$$ T^n(F)(x) = \frac{F(x)}{\lambda^n} - \sum\frac{b_i}{x -
c_i}\left(T^{n-1}(F)(R(c_i)) + \frac{T^{n - 2}(F)(R(c_i))}{\lambda} + ...
+ \frac{F(R(c_i))}{\lambda^{n - 1}}\right).$$
Continue calculation
$$ L_{m}(F) = \frac{1}{m}\sum_{n = 0}^{m - 1}T^n(F)(x) =
\frac{1}{m}\left(\sum_{n = 0}^{m-1}\frac{F(x)}{\lambda^n} - \sum\frac{b_i}{x
- c_i}\left(\sum_{n=0}^{m -1}\sum_{j=0}^n\frac{T^j(F)(d_i)}{\lambda^{n - j}}\right)\right).$$
For $ m\to\infty $ the first term tends to zero. Now let us calculate
the second term
$$\multline
\sum_{n=0}^{m -1}\sum_{j=0}^n\frac{T^j(F)(d_i)}
{\lambda^{n - j}}) = F(d_i)\left(\frac{1}{\lambda} + ...+ 
\frac{1}{\lambda^{m-1}}\right) +\\
+ T(F)(d_i)\left(\frac{1}{\lambda} + ... +
\frac{1}{\lambda^{m-2}}\right) + ... + T^{m-1}(F)(d_i)\frac{1}{\lambda} =
\endmultline
$$
$$ = F(d_i)\frac{\lambda}{\lambda - 1}\left(1 - \frac{1}{\lambda^{m-1}}\right) +
T(F)(d_i)\frac{\lambda}{\lambda - 1}\left(1 - \frac{1}{\lambda^{m-2}}\right) +
... + T^{m-1}(F)(d_i)\frac{1}{\lambda} = $$
$$ = \frac{\lambda}{\lambda - 1}\left(\sum_{j = 0}^{m-1} T^j(F)(d_i) -
\sum_{j=0}^{m-1}\frac{T^j(F)(d_i)}{\lambda^{m-i}}\right). $$
Such as $ \vert T^j(F)(d_i)\vert \leq M(F) $ for all $ j $ and $ i $ we
have that for $ m\to\infty $ the expression $
\frac{1}{m}\sum_i\left(\frac{b_i}{x-c_i}\sum_{j=0}^{m-1}
\frac{T^j(F)(d_i)}{\lambda^{m-i}}\right) $ tends to 0. So we conclude that $
L_j $ converges if and only if the functionals $ l^i_j $ converges. By
using assumptions we conclude that sequences $ L_{n_m} $ converges
pointwise and it complete claim.
\par By claim we know that for any fixed $ F \in W $ the sequence
functions $ A_{n_m}(T)(F)(a) $ converges on periodic points from Julia
set. We know that family functions $ A_{n_m}(T)(F)(a) $
forms bounded equicontinuous family functions. Now show that this
family has unique limit point. Indeed let $ F_1 $ and $ F_2 $ be two
different limit functions for our family, then by the claim $ F_1(x) =
F_2(x) $ for any repulsive periodic point and hence $ F_1 = F_2 $ on
Julia set. Now remember that functions $ F_i $ are holomorphic on Fatou
set we obtain $ F_1(z) = F_2(z) $ for all $ z \in {\Bbb C} $.
\par Now by using the fact $ \gamma_a = a\tau_1 -(a-1)\tau_0 + \tau_a $
we obtain that $ A_{n_m}(R^\ast)(\gamma_a) $converges weakly and hence
$ A_N(R^\ast)(\gamma_a) $ converges strongly for any $ a \in {\Bbb C}.$
This completes first step.
\par {\bf Second step} Here we prove that *-condition implies additional
conditions of the first step. Namely let us denote by $ Y $ the subset
of elements from $ L_1(J(R)) $ on which averages $ A_N(R^\ast) $ are
convergent. Note that $ Y $ is closed space such as family $ A_N(R^\ast)
$ forms equicontinuous family of operators.
\enddemo
\par We claim that {\it for any $ d_i = R(c_i) $ the elements
$\gamma_{d_i} $ belong to $ Y. $}
\demo{Proof of the claim}Otherwise it exists a continuous functional $ L
$ on $ L_1(J(R)) $ and $ i_0 $ so that $ L(\gamma_{i_0}) \neq 0 $ and
$ Y \subset ker(L).$ Note that $ L $ is invariant functional (i.e.
$L(R^\ast(f)) = L(f)$) such as for any $ f \in L_1(J(R)) $ the element $
f - R^\ast(f) $ belongs to $ Y. $ Let $ \nu \in L_\infty(J(R)) $ be the
element corresponding to $ L, $ then $ \nu $ is fixed point for Beltrami
operator $ B_R $ and hence the function $ F_\nu(a) =\iint\nu\tau $ is
fixed point for operator $ T $ i.e.
$$ 
\frac{F_\nu(R(a))}{R'(a)} - \sum\frac{b_iF_\nu(d_i)}{a - c_i} =
F_\nu(a).
$$
Let $ d $ be point from condition (*), then by lemma 22 and
under assumption $ F_\nu(x) = 0 $ for any $ x \in
\cup_{i=0}^\infty R^i(d) $. Therefore  meromorphic function $ \Phi(a) =
\sum\frac{b_iF_\nu(d_i)}{a - c_i} $ has big number ($> 2d-2$) of zeros
that immediately implies $\Phi(a) \equiv 0.$ Function $ F_\nu $ satisfy
to equation
$$ \frac{F_\nu(R(a))}{R'(a)} = F_\nu(a). $$
Finally we have that $ F_\nu $ is zero on set of all repulsive periodic
points hence on Julia set and hence everywhere because $ F_\nu $ is
holomorphic on Fatou set. Contradiction. Theorem is proved.\enddemo
\proclaim{Theorem 24} Let a map $ R $ be strongly convergent rational map.
Assume that Lebesgue measure of $PC(R)$ is zero. Then $ J(R) $ does not
support non-trivial invariant measurable fields.
\endproclaim
\demo{Proof} Theorem above and proposition 13 complete proof of 
theorem.
\enddemo
\proclaim{Corollary 25} Let rational map $ R $ be as in theorem above. 
Assume in addition that $ R $ is structurally stable, then
$ R $ is hyperbolic.
\endproclaim
\demo{Proof} Theorem $ A $ and Sullivan result (see
\cite{\bf MSS}) complete theorem.\enddemo
\par Further we will give sufficient conditions on rational maps to be strongly ergodic. This conditions will be given in terms of Poincare series of rational map. We start now with the following calculations.
\proclaim{Lemma 26} Let $ R $ be a rational map with no critical
relations and simple critical points. Let $ c $ be a critical point of $
R $ and $ d \in (R^k)^{-1}(c) $ be any point for some fixed $ k. $ Then
for any fixed $ m $ the coefficient $ b $ corresponding to the item $
\frac{1}{z - d} $ in expression $ s_m(z) $ has the following type
$$
b = \frac{1}{(R^m)''(d)} = \frac{1}{(R''(c))(R^{m - k -
1})'(R(c))((R^k)'(d))^2}.
$$
\endproclaim
\demo{Proof} Proof consists of consideration of residue in the point $
d. $ So let $ U $ be such neighborhood of $ c $ so that
\roster
\item  restriction $ R\vert_U : U \rightarrow R(U) $ presents 2-to-1
branching covering and
\item $ R(U) $ is the disc $\{z, \vert z - R(c)\vert = \epsilon\}$ for
some arbitrary fixed $ \epsilon. $
\item Let $ l \subset R(U) $ be an arc going from $ R(c) $ to $ \partial
R(U) $, then by $ B_1 $ and $ B_2 $ denote branches of $ R^{-1} $
mapping $ R(U) \backslash l $ into $ U $ and
\item in $ U $ the following decomposition is true
$ R(z) = R(c) + A(z - c)^2  + ...$
\endroster
\par Now let $ g $ be Jordan curve around point $ d $ eventually
mapping on $ \partial R(U). $ Write $ b =
\frac{1}{2i\pi}\int_g\frac{\partial z}{(R^m)'(z)}. $ Under conditions
there are no critical relations so it exists a branch $ J $ of $
(R^k)^{-1} $  such that $ J(c) = d$ and we have
$$ b = \frac{1}{2i\pi}\int_{\partial U}\frac{(J')^2(z)\partial z}{(R^{m -
k})(z)}                                              $$
and continue calculations
$$ b = \frac{1}{2i\pi}\int_{\partial R(U)}\frac{(J')^2(B_1^\prime(z))^2
+ (J')^2(B_2^\prime(z))^2}{(R^{m - k - 1})'(z)}\partial z. $$
Now remember that $ (B_1^\prime)^2(z) = \frac{1}{(R'(J(z)))^2} =
\frac{1}{(2A(B_1(z)- c) + ...)^2} = \frac{1}{4A^2(B_1(z) - c)^2 + ...}$
and use the fact
$$ 4A^2(B_1(z) -c)^2 = 4A(R(B_1(z)) - R(c) + ...) = 4A(z - R(c) + ...),
$$
and for $ z \in \partial R(U) $ the members $ ... $ are equivalent to $
\epsilon $. The same calculations for $ B_2 $ gives
$$ (B_2^\prime(z))^2 = \frac{1}{4A(z - R(c) + ...)}. $$
So we have
$$ b = \frac{1}{2i\pi}\int_{\partial R(U)}\frac{(J')^2(B_1(z))}{(4A(z -
R(c) + ...))(R^{m - k - 1})'(z)} + \frac{(J')^2(B_2(z))}{(4A(z -
R(c) + ...))(R^{m - k - 1})'(z)} $$
and now using arbitrariness of $ \epsilon $ obtain
$$ b = \frac{1}{2i\pi}\int_{\partial R(U)}\frac{(J')^2(B_1(z)) +
 (J')^2(B_2(z))}{(4A(z - R(c)))(R^{m - k - 1})'(z)}. $$

\par Further note that numerator under integral forms holomorphic
function in $ R(U) $ and $ (R^{m - k - 1})' $ does not have zeros in $
R(U) $ (otherwise we do $\epsilon$ smaller) and $ A = \frac{R''(c)}{2!}
$ we obtain
$$
b  = \frac{1}{(R''(c))(R^{m - k -
 1})'(R(c))((R^k)'(d))^2} = \frac{1}{(R^m)''(d)}.
$$
Lemma is proved
\enddemo
\par Follow by McMullen (\cite{\bf MM}) we recall the backward and forward Poincare series for the given rational map $ R.$
\proclaim{Definition} Forward Poincare series $ S(x, R) $
$$ P(x, R) = \sum_{n = 0}^\infty\frac{1}{\vert (R^n)'(R(x))\vert}. $$
Backward Poincare series $ P(x, R). $
Let $ \vert R^*\vert = R_{1,1}$ be the modulus of the Ruelle operator, 
then
$$ S(x, R) = \sum _{n = 1}^\infty \vert R^*\vert^n({\bold 1}_{\pmb{\Bbb C}})(x) = \sum_{n = 1}^\infty\sum_{R^n(y) = x}\frac{1}{\vert (R^n)'(y)\vert^2}.$$
\endproclaim

Let us again consider the function $ s_n(a)=\sum\frac{\vert b_i\vert}{\vert a - c_i\vert} $ and let  $ B_n = \sum\vert b_i\vert. $ Then by lemma above we have.
$$ B_n = \sum_{c \in Cr(R)}\frac{1}{\vert R''(c)\vert}\sum_{j = 1}^{n -1}\frac{1}{\vert (R^{n - j - 1})'(R(c))\vert}\sum_{R^j(y) = c}\frac{1}{\vert(R^j)'(y)\vert^2} $$
and hence we have the following formal equality
$$ \sum_{n = 2} B_n = \sum_{c \in Cr(R)}\frac{1}{\vert R''(c)\vert}S(c, R)\otimes P(c, R),
$$
we recall that $\otimes $ means Cauchy product of series.
\proclaim{Corollary 27} Let $ R $ be a rational map with simple critical points and no simple critical relation. Then for  any $ a \in F(R) $ the function $ s_n(a) $ is bounded in the following cases.
\roster
\item Collet-Eckmann case. For any critical point $ c \neq \infty $ the series $ S(c,R) $ and $ P(c, R) $ are bounded.
\item For any critical point $ c \neq \infty  $ one of the series  $ S(c,R) $ or $ P(c, R) $ are bounded and the second one has uniformly bounded elements.
\item Conjectural case. Both series diverge slow enough (like harmonic series).
\endroster
\endproclaim
\demo{Prof} The cases (1) and (2) are immediately follows from properties of Cauchy product. In last case also follows from properties of Cauchy product, such as Cauchy product of two harmonic series is divergent but has uniformly bounded elements. Let us again repeat that it is not clear does exists a rational map for which both Poincare series are equivalent to harmonic series.
\enddemo
\par Let us note that this corollary looks like corollary A with Poincare series.
\subheading{Convergent maps}Now the our aim is to give a weaker condition on a rational map to be ergodic.
	Let us recall that the main goal of the discussion above is the estimating the norm of the operator $ T $ on space $ X, $ that is to estimate the expression
	$$ \sum\frac{b_iF(R(c_i))}{a - c_i}. $$
Let us rewrite this expression by using lemma 26.
\par Let $ R $ be a rational map and $ c_i, i =1,..., 2deg(R) -2, $ and $ d_i, i =1,...,  2deg(R) -2, $ be critical points and critical values, respectively and let point $ z  =\infty $ be superattractive point. Then by induction we have.
$$ \frac{1}{R'(z)} = \sum_i\frac{b_i}{z - c_i} = \sum_i\frac{1}{R''(c_i)}\frac{1}{z - c_i}, $$
$$ ... $$
$$ \frac{1}{(R^n)'(z)} = \sum_i\sum_{k = 0}^{n - 1}\left(\sum_{y \in R^{-k}(c_i}\frac{1}{(R^{n})''(y)}\frac{1}{z - y}\right). $$
Hence we obtain
$$ 
\sum\frac{b^n_iF(R(c_i))}{a - c^n_i} = \sum_i\sum_{k = 0}^{n - 1}F(R^{n - k - 1}(d_i))\left(\sum_{y \in R^{-k}(c_i)}\frac{1}{(R^{n})''(y)}\frac{1}{z - y}\right). 
$$
\proclaim{Lemma 28} 
$$ 
\multline
\sum_{y \in R^{-k}}\frac{1}{(R^{n})''(y)}\frac{1}{a - y} = 
\frac{1}{R''(c_i)}\frac{1}{(R^{n - k - 1})'(d_i)}\sum_j\frac{(J_j')^2(c_i)}{a - J_j(c_i)} =\\ =\frac{1}{R''(c_i)}\frac{1}{(R^{n - k - 1})'(d_i)}(R^*)^k(-\tau_a)(c_i),
\endmultline
$$
where $ J_j $ are branches of $ R^{-k}, \tau_a(z) = \frac{1}{z - a} $ and $ R^* $ is Ruelle operator.
\endproclaim
\demo{Proof} Lemma 26 and equalities above complete this lemma.
\enddemo
Then we have the proposition.
\proclaim{Proposition A} Let $ R $ be rational map with simple critical points. Let $\infty $ be a fixed point for $ R. $ Then there exist the following formal relations. 
$$ \align RP(a, R) - 1
&= \sum_i \lambda ^i -\sum_i\frac{1}{R''(c_i)}RS(c_i,R,a)\otimes RP(c_i,R),\text{ where } \lambda \text{ is multiplier of } \infty\\
RS(x, R, a) 
&= A(x, R, a) - \sum_k\frac{1}{R''(c_k)}A(c_k,R, a)\otimes RS(x,R,c_k),
\endalign
$$
where $ c_k $ are critical points of $ R. $
\endproclaim
\demo{Proof} The first equality is immediate corollary of discussion above. 
\par Let us show the second equality. By lemma 14 we can calculate.
$$ 
\multline
\left(R^*\right)^0(\tau_a)(z) =\tau_a(z), \left(R^*\right)(\tau_a)(z) = \frac{1}{R'(a)(z - R(a)} - \sum_i\frac{b_i}{(a - c_i)(z - R(c_i))}\\
\left(R^*\right)^2(\tau_a)(z) = \frac{1}{(R^2)'(a)(z - R^2(a)} - \frac{1}{R'(a)}\sum_i\frac{b_i}{(R(a) - c_i)(z - R(c_i)} - \\
\sum_i\frac{b_i}{(a - c_i)}R^*(\frac{1}{(z - R(c_i))}
\endmultline
$$
and by induction
$$\multline
\left((R^*)\right)^n(\tau_a)(z) = \frac{1}{(R^n)'(a)(z - R^n(a)} - 
\sum_ib_i\left(\frac{1}{(R^{n -1})'(a)(R^{n - 1}(a) - c_i)(z - R(c_i)}\right). +\\
  ... + \left(. \frac{1}{a - c_i}(R^{n - 1})^*(\frac{1}{z - R(c_i)}\right)
\endmultline
$$
hence summation respect to $n $ gives desired equality. Lemma is proved.
\enddemo
 Concluding the discussion above we obtain the following expression.
$$ 
\frac{1}{(R^n)'(a)} = -\sum_i\sum_{k = 0}^{n - 1}\frac{1}{R''(c_i)}\frac{1}{(R^{n - k - 1})'(d_i)}(R^*)^k(\tau_a)(c_i)
$$ 
and
$$
\sum\frac{b^n_iF(R^n(c_i))}{a - c_i} = -\sum_i\sum_{k =0}^{n -1}\frac{F(R^n(c_i))}{R''(c_i)(R^{n - k - 1})'(c_i)}(R^*)^k(\tau_a)(c_i).
$$

\proclaim{Theorem 29} Let $ R $ be convergent map. Then $ R $ is mean ergodic.
\endproclaim
\demo{Proof} Assume $ R $ satisfies itself to condition $ * $ from definition of convergent map. Then by using arguments of  theorem 23 and lemma 22 we have
$$ \vert \frac{1}{N_i}\sum_{j = 0}^{N_i} T^j(F)(x)\vert \leq 2 A_{N_i}(x, R,)\Vert F\Vert_\infty.
$$
Hence sequences $ \frac{1}{N_i}\sum_{j = 0}^{N_i} T^j(F)(x) $ is convergent (up to passing to a subsequences) in $*-$weak topology on the space $ W $ for any $ x \in \cup_nR^n(a). $ Again by using arguments  theorem 23 and Mean ergodicity lemma we conclude that $ R $ is mean ergodic.
\enddemo
\proclaim{Theorem A} Let $ R $ be convergent map. Assume Lebesgue measure of postcritical set is zero, then there is no non-trivial invariant line field on Julia set.
\endproclaim
\demo{Proof} Theorem 29 and proposition 13 give desired conclusion.
\enddemo
In corollary 27 we use definition of Collet-Eckmann maps in Przytycki sense. That is for all critical points the both forward and backward Poincare series are absolutely convergent. Now we redefine these maps by the following way.
\proclaim{Definition} A rational map $ R $ we will call {\rm Collet-Eckmann} map if does exists a point $ a \in {\Bbb C} $ with long orbit,
$ \#\{\cup_nR^n(a)\} > 2deg(R) - 1 $ such that\hfill\newline
{\rm For any critical point $ c $ and any $ x \in  \{\cup_nR^n(a)\} $ the both Ruelle-Poincare series $ RS(c, R, x) $ and $ RP(c, R) $ are absolutely convergent.}
\endproclaim
\proclaim{Corollary A}
Let $ R $ be a rational map with simple critical points and no simple critical relation. Then $ R $ is convergent map if for some $ a \in {\Bbb C} $ with $\#\{\cup_iR^i(a)\} > 2deg(R) - 1 $ the one of the following is true.
\roster
\item Collet-Eckmann maps.
\item For any critical point $ c  $ and $ x \in \{\cup_iR^i(a)\}  $ one of the series  $ RS(c,R, x) $ and $ RP(c, R) $ is  absolutely convergent and second one has uniformly bounded elements.
\item Conjectural case. Both series diverge slow enough (like harmonic series).
\endroster
\endproclaim
\demo{Proof} Evidently the maps from all (1) - (3) cases are convergent.
\enddemo 


\heading{\bf Measures}
\endheading
\par Start again with  rational map $ R. $ 
Consider an element $ \gamma \in A(R) $ and Cesaro average sequence $ A_N(R)(\gamma)$ $=$ $\frac{1}{N}\sum_{i = 0}^{N - 1} (R^*)^i(\gamma). $ Let $ C(U) $ be the space of continuous functions defined on $ \overline U $ for the fixed essential neighborhood $ U. $ Then any $*$-weak limit of $ A_N(R)(\gamma) $ on $ C(U) $ we call {\it weak boundary} of $  \gamma $ respect to $ R^* $ over $ U $ and denote set of all limit measures  by $ \gamma(U,R).$ 
\proclaim{Theorem 30} Let $ R $ be a structurally stable rational map with non empty Fatou set. Assume there exists non-zero weak boundary $ \mu \in \gamma(U, R^*) $ for an element $ \gamma \in A(R) $ and an essential neighborhood $ U. $ Then Lebesgue measure
$ m(J(R)) > 0 $ and there exists a non-trivial invariant line field on $ J(R).$
\endproclaim
\demo{Proof} Under assumptions there exists an essential $ U $ and $ \gamma \in A(R) $ and  subsequence $ N_i $ such that 
\roster
\item $ \iint\phi A_{N_i}(R)(\gamma) $ converges for any $ \phi \in C(U) $ and
\item there exists $ \psi \in C(U) $ such that $\lim_{i\rightarrow\infty} \iint\psi A_{N_i}(R)(\gamma) \neq 0. $
\endroster
	By using density of space of compactly supported continuous function at the space $ C(U) $ we can assume that $ \psi $ has a compact support $ D \subset \overline{U}. $ Continue $ \psi $ on $ \overline{\Bbb C} \backslash D $ by zero
we obtain $ \lim_{i\rightarrow\infty}\iint_{\overline{\Bbb C}}\psi A_{N_i}(R)(\gamma) \neq 0. $  Hence the dual average $ A_N(B_{R})(\psi) = \frac{1}{N}\sum_{i = 0}^{N - 1} (B_{R})^i(\psi) $ has non-zero $*$-weak limit element in $ *$-weak topology on $ L_\infty(J(R)). $ Let $ \mu \in L_\infty(J(R)) $ be this non-zero limit element, then $ \mu $ is fixed for $ B_{R} $ and $ \mu = 0 $ on $ F(R) $ by construction. Hence
$ m(J(R)) > 0 $ and $ \mu $ defines desired invariant line field.
\enddemo

It is not clear when the inverse statement is true. But we suggest the following conjecture.
\proclaim{Conjecture} Let $ R $ be a rational map with non-empty Fatou set, the $ T(J(R)) = \emptyset $ if and only if  weak boundaries $ \gamma(U,R^*) = 0 $ for all $ \gamma \in A(R) $ and all essential neighborhood $ U. $
\endproclaim
\par In general the absence of invariant line fields on Julia set means the mean ergodicity of $ R^* $ on $ L_1(J(R)) $ and so it is interesting to understand the conditions implying the mean ergodicity of $ R^* $ from measure's point of view. To do this let us recall definition of the following objects:
\roster
\item $ U $ is an essential neighborhood of $ J(R) $ and
\item $ H(U) $ consists of $ h\in C({\overline U}) $ such that $ \frac{\partial h}{\partial \overline{z}} $ (in sense of distributions) belongs to $ L_\infty(U) $
\item $ H(U) $ inherits the topology of $ C({\overline U}).$
\endroster
\par Measures $ \nu^i_l. $
\roster
\item Let $ c_i $ and $ d_i $ be critical points and critical values, respectively.
Then  define $ \mu_n^i = \frac{\partial}{\partial\overline{z}}((R^*)^n(\frac{1}{z - d_i})) $ (in sense of distributions). We will show below that $  (R^*)^n(\frac{1}{z - d_i})$ $ = $ $\sum_{j = 0}^n\frac{\alpha^i_j}{z - R^j(d_i)} $ and hence $ \mu_n^i = \sum_{j = 0}^n\alpha^i_j\delta_{R^j(d_i)}, $ where $ \delta_a $ denotes the delta measure with mass at the point $ a. $
\item Define by $ \nu^i_l $ the average $ \frac{1}{l}\sum_{k = 0}^{l - 1}\mu_k^i. $
\endroster

\proclaim{Theorem C} Let $ R $ be a rational map with simple critical points and no critical relations. Assume that $ F(R) \neq \emptyset $ and $ m(Pc(R)) = 0, $ where $ Pc(R)$ is the postcritical set and $ m $ is denote the Lebesgue measure. Then $ T(J(R)) = \emptyset $ if and only if there exist an essential neighborhood $ U $ and a sequences of integers $ \{l_k\} $ such that the measures $ \{\nu^i_{l_k}\} $ converges in $*$-weak topology on $ H(U) $ for any $ i = 1, ..., 2deg(R) - 2.$
\endproclaim
\demo{Proof} Suppose that $ \nu^i_{l_k} $ converges  in $*$-weak topology on $ H(U) $ for all $ i $  a subsequence $ \{l_k\} $ and an essential neighborhood $ U. $ Then the sequence averages $ A_N(R)(\frac{1}{z - d_i}) \in L_1(U) $ is convergent weakly, in case $ m(J(R)) > 0 $ that means  $ A_N(R)(\frac{1}{z - d_i}) $ converges strongly in $ L_1(J(R)). $ By using arguments proposition 13 and theorem 23 (second step) there exist no invariant line fields.
\par Now assume there exist no invariant line fields on  $ J(R) $ Let us show that $ \nu_l^i \rightarrow 0 $ in $ *$-weak topology on $ H(U) $ for any essential neighborhood $ U. $ Otherwise there exist a sequence $ \{l_k\} $ an essential neighborhood $ U $ and a  function $ F \in H(U) $ such that 
$$ \lim_{k\rightarrow\infty}\iint F\nu^{i_0}_{l_k} = \lim_{k\rightarrow\infty}\iint F_{\overline{z}}\frac{1}{l_k}\sum_{n = 0}^{l_k - 1} (R^*)^n\left(\frac{1}{z - d_{i_0}}\right) \neq 0. $$
Again like in previous theorem we can think that $ F $ is continuous differentiable function with compact support belonging to $ U. $ Hence the sequences $ \frac{1}{l_k}\sum_{i = 0}^{l_k - 1}B_R(F_{\overline{z}}) $ has non-zero limit element $ \mu \in L_\infty(U) $ in $*$-weak topology on $ L_\infty(U). $ By arguments of Mean ergodicity lemma $ \mu $ is invariant. Contradiction. Lemma is proved. 
\enddemo

\enddocument